\documentclass[10pt]{amsart}
 
\usepackage{amscd,a4,graphics,amssymb,amsbsy}
\usepackage[all]{xy}
\usepackage{color}

\newtheorem{thm}{Theorem}[section] \newtheorem{lemma}[thm]{Lemma}
\newtheorem{prop}[thm]{Proposition}
\newtheorem{cor}[thm]{Corollary}

\theoremstyle{definition}
 \newtheorem{dfn}[thm]{Definition}
\newtheorem{ntn}[thm]{Notation} \newtheorem{rmk}[thm]{Remark}
\newenvironment{pf}{\medskip\noindent{{\em Proof: }}}{\qed}
\newtheorem{ex}[thm]{Example}
\newtheorem{ass}[thm]{Assumption}

\newcommand{\Spec}{\mathrm{Spec}}

\newcommand{\Char}{\mathrm{char}}

\newcommand{\Gal}{\mathrm{Gal}}

\newcommand{\Gcd}{\mathrm{gcd}}
\newcommand{\Deg}{\mathrm{deg}}
\newcommand{\Mult}{\mathrm{mult}}
\newcommand{\Lcm}{\mathrm{lcm}}

\newcommand{\Sing}{\mathrm{Sing}}
\newcommand{\Aut}{\mathrm{Aut}}
\newcommand{\Id}{\mathrm{id}}

\newcommand{\Dim}{\mathrm{dim}}
\newcommand{\Ker}{\mathrm{Ker}}

\newcommand{\Card}{\mathrm{Card}}
\newcommand{\Pic}{\mathrm{Pic}}

\newcommand{\Tr}{\mathrm{Tr}}

\numberwithin{equation}{section}

\begin{document}


\title{Galois actions on N\'eron models of Jacobians}


  \author{Lars Halvard Halle}\email{larshalvard.halle@wis.kuleuven.be}
  \address{Departement Wiskunde, Celestijnenlaan 200B, B-3001 Leuven (Heverlee), Belgium}
    

\begin{abstract}
Let $X$ be a smooth curve defined over the fraction field $K$ of a complete d.v.r. $R$, and let $K'/K$ be a tame extension. We study extensions of the $G = \Gal(K'/K)$-action on $ X_{K'} $ to certain regular models of $X_{K'}$ over $R'$, the integral closure of $R$ in $K'$. In particular, we consider the induced action on the cohomology groups of the structure sheaf of the special fiber of such a regular model, and obtain a formula for the Brauer trace of the endomorphism induced by a group element on the alternating sum of the cohomology groups.  

We apply these results to study a natural filtration of the special fiber of the N\'eron model of the Jacobian of $X$ by closed, unipotent subgroup schemes. We show that the jumps in this filtration only depend on the fiber type of the special fiber of the minimal regular model with strict normal crossings for $X$ over $R$, and in particular are independent of the residue characteristic. Furthermore, we obtain information about where these jumps occur. We also compute the jumps for each of the finitely many possible fiber types for curves of genus $1$ and $2$.   
\end{abstract}

\thanks{The research was partially supported by the Fund for Scientific Research - Flanders (G.0318.06)}
\maketitle

\section{Introduction}

Let $X$ be a smooth, projective and geometrically connected curve of genus $ g(X) > 0 $, defined over the fraction field $K$ of a complete discrete valuation ring $R$, with algebraically closed residue field $k$. By a \emph{model} for $X$ over $R$, we mean an integral and normal scheme $ \mathcal{X} $ that is flat and projective over $ S = \Spec(R) $, and with generic fiber $ \mathcal{X}_K \cong X $. The special fiber $ \mathcal{X}_k $ of such a model is called a \emph{reduction} of $X$. 

The semi-stable reduction theorem, due to Deligne and Mumford (\cite{DelMum}, Corollary 2.7), states that there exists a finite, separable field extension $ L/K$ such that $X_L$ admits a semi-stable model over the integral closure $R_L$ of $R$ in $L$.

In order to study reduction properties of $X$, it can often be useful to work with the Jacobian $J/K$ of $X$. The question whether $X$ has semi-stable reduction over $S = \Spec(R)$ is reflected in the structure of the \emph{N\'eron model} $ \mathcal{J}/S$ (cf. ~\cite{Ner}) of $J$. In fact, $X$ has semi-stable reduction over $S$ if and only if $ \mathcal{J}_k^0 $, the identity component of the special fiber, has zero \emph{unipotent radical} (\cite{DelMum}, Proposition 2.3).

In general, it is necessary to make ramified base extensions in order for $X$ to obtain semi-stable reduction. If the residue characteristic is positive, it can often be difficult to find explicit extensions over which $X$ obtains stable reduction. In the case where a tamely ramified extension suffices one can do this by considering the geometry of suitable regular models for $X$ over $ S $ (cf. ~\cite{Thesis}). In this paper we study, among other things, how the geometry of the N\'eron model contains information that is relevant for obtaining semi-stable reduction for $X$.

\subsection{N\'eron models and tame base change}
Let $ K'/K $ be a finite, separable and tamely ramified extension of fields, and let $R'$ be the integral closure of $R$ in $ K' $. Then $R'$ is a complete discrete valuation ring, with residue field $k$. Furthermore, $ K'/K $ is Galois, with group $ G = \boldsymbol{\mu}_n $, where $ n = \Deg(K'/K) $.

Let $ \mathcal{J}'/S' $ be the N\'eron model of the Jacobian of $ X_{K'} $, where $ S' = \Spec(R') $. Due to a result by B. Edixhoven (\cite{Edix}, Theorem 4.2), it is possible to describe $ \mathcal{J}/S $ in terms of $ \mathcal{J}'/S' $, together with the induced $G$-action on $ \mathcal{J}' $. Namely, if $W$ denotes the \emph{Weil restriction} of $ \mathcal{J}'/S' $ to $S$ (cf. ~\cite{Ner}, Chapter 7), one can let $G$ act on $W$ in such a way that $ \mathcal{J} \cong W^G $, where $ W^G $ denotes the scheme of invariant points. In particular, one gets an isomorphism $ \mathcal{J}_k \cong W_k^G $. By \cite{Edix}, Theorem 5.3, one can use this description of $ \mathcal{J}_k $ to define a descending filtration 
$$ \mathcal{J}_k = F_n^0 \supseteq \ldots \supseteq F_n^i \supseteq \ldots \supseteq F_n^n = 0 $$
of $ \mathcal{J}_k $ by closed subgroup schemes.

In \cite{Edix}, Remark 5.4.5, a generalization of this setup is suggested. If we define $ \mathcal{F}^{i/n} = F_n^i $, where $ F_n^i $ is the $i$-th step in the filtration induced by the extension of degree $n$, one can consider the filtration 
$$ \mathcal{J}_k = \mathcal{F}^0 \supseteq \ldots \supseteq \mathcal{F}^a \supseteq \ldots \supseteq \mathcal{F}^1 = 0, $$
with indices in $ \mathbb{Z}_{(p)} \cap [0,1] $. The construction of $ \mathcal{F}^a $ is independent of the choice of representatives $i$ and $n$ for $ a = i/n $.

The filtration $ \{ \mathcal{F}^a \} $ contains significant information about $\mathcal{J}$. For instance, the subgroup schemes $ \mathcal{F}^a $ are \emph{unipotent} for $ a > 0 $, so in a natural way, this filtration gives a measure on how far $ \mathcal{J}/S $ is from being \emph{semi-abelian}.

One way to study the filtration $ \{ \mathcal{F}^a \} $ is to determine where it \emph{jumps}. This will occupy a considerable part of this paper. The jumps in the filtration often give explicit numerical information about $X$. For instance, if $X$ obtains stable reduction after a tamely ramified extension, we show that the jumps occur at indices of the form $i/\tilde{n}$, where $\tilde{n}$ is the degree of the minimal extension that realizes stable reduction for $X$. 

It follows from Edixhoven's theory that to determine the jumps in the filtration $ \{ F^i_n \} $ induced by an extension of degree $n$, one needs to compute the irreducible characters for the representation of $ \boldsymbol{\mu}_n $ on the tangent space $ T_{ \mathcal{J}'_k, 0 } $. We shall use such computations for \emph{infinitely} many integers $n$ to describe the jumps of the filtration $ \{ \mathcal{F}^a \} $ with rational indices.

\subsection{N\'eron models for Jacobians}
Contrary to the case of general abelian varieties, N\'eron models for Jacobians can be constructed in a fairly concrete way, using the theory of the relative Picard functor (cf. ~\cite{Ner}, Chapter 9). The following property will be of particular importance to us: If $ \mathcal{Z}/S' $ is a regular model for $ X_{K'}/K' $, then there is a canonical isomorphism
$$ \Pic_{\mathcal{Z}/S'}^0 \cong (\mathcal{J}')^0, $$  
where $ \Pic_{\mathcal{Z}/S'}^0 $ (resp.~$(\mathcal{J}')^0$) is the identity component of $ \Pic_{\mathcal{Z}/S'} $ (resp.~$\mathcal{J}'$). It follows that there is a canonical isomorphism
$$ H^1(\mathcal{Z}_k, \mathcal{O}_{\mathcal{Z}_k}) \cong T_{\mathcal{J}'_k, 0}. $$

We shall work with regular models $ \mathcal{Z} $ of $ X_{K'} $ that admit $G$-actions that are compatible with the $G$-action on $ \mathcal{J}' $. It will then follow that the representation of $G$ on $ T_{\mathcal{J}'_k, 0} $ can be described in terms of the representation of $G$ on $ H^1(\mathcal{Z}_k, \mathcal{O}_{\mathcal{Z}_k}) $.

\subsection{Models and actions}
In order to find an $S'$-model for $X_{K'} $ with a compatible $G$-action, we take a model $ \mathcal{X} $ of $X$ over $S$, and consider its pullback $ \mathcal{X}_{S'} $ to $ S' $. Let $ \mathcal{Y} \to \mathcal{X}' \to \mathcal{X}_{S'} $ be the composition of the normalization with the minimal desingularization. Then $ \mathcal{Y} $ is a model of $ X_{K'} $ with an action of $G$ that lifts the obvious action on $ \mathcal{X}_{S'} $. The $G$-action restricts to the special fiber $ \mathcal{Y}_k $, and in particular, $G$ will act on the cohomology groups $ H^i(\mathcal{Y}_k, \mathcal{O}_{\mathcal{Y}_k}) $, for $ i = 0, 1 $. 

In order to understand the $G$-action on $ H^i(\mathcal{Y}_k, \mathcal{O}_{\mathcal{Y}_k}) $, we need a good description of the geometry of $\mathcal{Y}$ and of the $G$-action on $ \mathcal{Y} $. For this purpose, we demand that the model $ \mathcal{X} $ has good properties. To begin with, we shall require that $ \mathcal{X} $ is regular, and that the special fiber is a divisor with strict normal crossings. Furthermore, we shall always require that any two irreducible components of $ \mathcal{X}_k $, whose multiplicities are both divisible by the residue characteristic, have empty intersection. This condition is automatically fulfilled if $X$ obtains stable reduction after a tamely ramified extension, but holds also for a larger class of curves.

Under these assumptions, it turns out that the normalization $ \mathcal{X}' $ of $ \mathcal{X}_{S'} $ has at most \emph{tame cyclic quotient singularities} (cf. ~\cite{CED}, Definition 2.3.6 and \cite{Thesis}, Paper I Proposition 4.3). These singularities can be resolved explicitly, and it can be seen that $ \mathcal{Y} $ is a strict normal crossings model for $ X_{K'} $.

We shall also only consider the case where $ n = \Deg(K'/K) $ is relatively prime to the multiplicities of all the irreducible components of $ \mathcal{X}_k $. With this additional hypothesis, it turns out that we can describe the combinatorial structure of the special fiber $ \mathcal{Y}_k $ (i.e., the intersection graph of the irreducible components, their genera and multiplicities), in terms of the corresponding data for $ \mathcal{X}_k $.

If all the assumptions above are satisfied, it follows that all irreducible components of $ \mathcal{Y}_k $ are stable under the $G$-action on $ \mathcal{Y} $, and that all intersection points in $ \mathcal{Y}_k $ are fixed points. We can explicitly describe the action on the cotangent space of $ \mathcal{Y} $ at these intersection points, and the restriction of the $G$-action to each irreducible component of $ \mathcal{Y}_k $. 

\subsection{Action on cohomology}
Next, we study the representation of $ G = \boldsymbol{\mu}_n $ on $ H^1(\mathcal{Y}_k, \mathcal{O}_{\mathcal{Y}_k}) $. In particular, we would like to compute the irreducible characters for this representation. So for every $ g \in G $, we want to compute the trace of the endomorphism of $ H^1(\mathcal{Y}_k, \mathcal{O}_{\mathcal{Y}_k}) $ induced by $g$, and then use this information to find the characters.

There are some technical problems that need to be overcome in order to do this. First, since we allow the residue characteristic to be positive, just knowing the trace for each $g \in G$ may not give sufficient information to compute the characters. Instead, we have to compute the so called \emph{Brauer trace} for every $g \in G $ (cf. \cite{SerreLin}, Chapter 18). This means that we have to lift the eigenvalues and traces from characteristic $p$ to characteristic $0$. From knowing the Brauer trace for every $ g \in G $ we can compute the irreducible Brauer characters, and then the ordinary characters are obtained by reducing the Brauer characters modulo $p$. Second, the special fiber $ \mathcal{Y}_k $ will in general be singular, and even non-reduced. This complicates trace computations considerably.

To deal with these problems, we introduce in Section \ref{section 6} a certain filtration of the special fiber $ \mathcal{Y}_k $ by effective subdivisors, where the difference at the $i$-th step is an irreducible component $ C_i $ of $ \mathcal{Y}_k $. Since $ \mathcal{Y} $ is an SNC-model, each $ C_i $ is a smooth and projective curve, and with our assumption on $n$, the $G$-action restricts to each $ C_i $. Furthermore, to each step in this filtration, one can in a natural way associate an invertible $ G $-sheaf $ \mathcal{L}_i $, supported on $C_i$. 

We apply the so called Lefschetz-Riemann-Roch formula (\cite{Don}, Corollary 5.5), in order to get a formula for the Brauer trace of the endomorphism induced by each $ g \in G $ on the formal difference $ H^0(C_i,\mathcal{L}_i) - H^1(C_i,\mathcal{L}_i) $. An important step is to show that our description of the action on $ \mathcal{Y} $ is precisely the data that is needed to obtain these formulas. Then we show that these traces add up to give the Brauer trace for the endomorphism induced by each $ g \in G $ on the formal difference $ H^0(\mathcal{Y}_k, \mathcal{O}_{\mathcal{Y}_k}) - H^1(\mathcal{Y}_k, \mathcal{O}_{\mathcal{Y}_k}) $. In particular, we give in Theorem \ref{thm. 9.13}, one of the main results in this work, a formula for this Brauer trace, and show that it only depends on the combinatorial structure of $ \mathcal{X}_k $. 

Let us also remark that in our situation, we already know the character for $ H^0(\mathcal{Y}_k, \mathcal{O}_{\mathcal{Y}_k}) $, and hence we will be able to compute the irreducible characters for $ H^1(\mathcal{Y}_k, \mathcal{O}_{\mathcal{Y}_k}) $ in this way.

\subsection{Conclusions and computations}
If now $ \mathcal{X}/S $ is the minimal regular model with strict normal crossings for $ X/K $, we prove in Theorem \ref{main character theorem} that the irreducible characters for the representation of $ G = \boldsymbol{\mu}_n $ on $ H^1(\mathcal{Y}_k, \mathcal{O}_{\mathcal{Y}_k}) $ only depend on the combinatorial structure of the special fiber $ \mathcal{X}_k $, as long as $n$ is relatively prime to $l$, where $l$ is the least common multiple of the multiplicities of the irreducible components of $ \mathcal{X}_k $. 

Let $ \mathcal{J} $ be the N\'eron model of the Jacobian of $X$. Then it follows from Theorem \ref{main character theorem} that the jumps in the filtration $ \{ \mathcal{F}^a \mathcal{J}_k \} $ only depend on the combinatorial structure of $ \mathcal{X}_k $ (Corollary \ref{main jump corollary}). This is due to the fact that $ \mathbb{Z}_{ (p l) } \cap [0,1] $ is ``dense'' in $ \mathbb{Z}_{ (p) } \cap [0,1] $. Furthermore, in Corollary \ref{specific jump corollary}, we draw the conclusion that the jumps are actually independent of ~$p$, and that the jumps can only occur at finitely many rational numbers of a certain kind, depending on the combinatorial structure of $ \mathcal{X}_k $. 

For a fixed genus $ g \geq 1 $, there are only a finite number of possible combinatorial structures for $ \mathcal{X}_k $, modulo a certain equivalence relation. In case $ g = 1 $ or $ g = 2 $, one has complete classifications (cf. ~\cite{Kod} for $g=1$ and \cite{Ueno}, \cite{Ogg} for $g=2$). In Section \ref{computations and jumps} we compute the jumps for each possible fiber type for $g=1$ (which were also computed by Schoof in \cite{Edix}) and for $g=2$.

\subsection{Acknowledgements}
I would like to thank Bas Edixhoven for suggesting this subject to me, and generously sharing his ideas. I would also like to thank my thesis advisor Carel Faber for discussing the material in this paper with me.

\section{N\'eron models and tamely ramified extensions}\label{Neron}

\subsection{} 
Let $ R $ be a discrete valuation ring, with fraction field $K$ and residue field $k$, and let $A$ be an abelian variety over $K$. There exists a canonical extension of $A$ to a smooth group scheme $ \mathcal{A} $ over $ S = \Spec(R) $, known as the \emph{N\'eron model} (\cite{Ner}, Theorem 1.4/3). The N\'eron model is characterized by the following universal property: for every smooth morphism $ T \rightarrow S $, the induced map $ \mathcal{A}(T) \rightarrow A(T_K) $ is \emph{bijective}.

\subsection{}\label{neronbase}
We assume from now on that $R$ is strictly henselian. Let $ K'/K $ be a finite, separable extension of fields, and let $R'$ be the integral closure of $R$ in $K'$. Let $ \mathcal{A}'/S' $ denote the N\'eron model of the abelian variety $ A_{K'}/K' $, where $ S' = \Spec(R') $. In general, it is not so easy to describe how N\'eron models change under ramified base extensions. However, in the case where $K'/K$ is tamely ramified, one can relate $ \mathcal{A}'/S' $ and $ \mathcal{A}/S $ in a nice way, due to a result by B. Edixhoven (\cite{Edix}, Theorem 4.2). We will in this section explain this relation, following the treatment in \cite{Edix}. We refer to this paper for further details. 

Assume now that $ K'/K $ is tamely ramified. Then $ K'/K $ is Galois with group $ G = \boldsymbol{\mu}_n $, where $ n = [K' : K] $. Let $G$ act on $ A_{K'} = A \times_{\Spec(K)} \Spec(K') $ (from the right), via the action on the right factor. By the universal property of $ \mathcal{A}' $, this $G$-action on $ A_{K'} $ extends uniquely to a right action on $ \mathcal{A}' $, such that the morphism $ \mathcal{A}' \rightarrow S' $ is equivariant. The idea in \cite{Edix} is to reconstruct $ \mathcal{A} $ as an invariant scheme for this action. However, since $ \mathcal{A} $ is an $S$-scheme, it is necessary to ``push forward'' from $S'$ to $S$. 

\subsection{}

The \emph{Weil restriction} of $ \mathcal{A}' $ to $S$ is the functor $ \Pi_{S'/S}(\mathcal{A}'/S') $ defined by assigning, to any $S$-scheme $T$, the set $ \Pi_{S'/S}(\mathcal{A}'/S')(T) = \mathcal{A}'(T_{S'}) $ (cf. ~\cite{Ner}, Chapter 7). This functor is representable by an $S$-scheme, which we will denote by $ \Pi_{S'/S}(\mathcal{A}'/S') $.


In \cite{Edix}, an equivariant $G$-action on $ \Pi_{S'/S}(\mathcal{A}'/S') $ is defined, corresponding to the $G$-action on $ \mathcal{A}' $.
Furthermore, there is a canonical morphism
$$ \mathcal{A} \rightarrow \Pi_{S'/S}(\mathcal{A}'/S'), $$
which, according to \cite{Edix}, Theorem 4.2, is a closed immersion, and induces an isomorphism
\begin{equation}\label{Bastheorem}
\mathcal{A} \cong (\Pi_{S'/S}(\mathcal{A}'/S'))^G, 
\end{equation}
where $(\Pi_{S'/S}(\mathcal{A}'/S'))^G$ denotes the \emph{scheme of invariant points} for this $ G $-action.



\subsection{Filtration of $\mathcal{A}_k$} 
One can use Isomorphism (\ref{Bastheorem}) to study the special fiber $ \mathcal{A}_k $ in terms of $ \mathcal{A}'_k $, together with the $G$-action. Indeed, let $ R \subset R' = R[\pi']/(\pi'^n - \pi) $ be a tame extension, where $ \pi $ is a uniformizing parameter for $R$. Then we have that $ R'/\pi R' = k[\pi']/(\pi'^n) $. For any $ k $-algebra $C$, it follows that 
$$ \mathcal{A}_k(C) \cong X_k^G(C) \cong X_k(C)^G \cong \mathcal{A}'(C[\pi']/(\pi'^n))^G, $$
where $ X = \Pi_{S'/S}(\mathcal{A}'/S') $.

In \cite{Edix}, Chapter 5, this observation is used to construct a filtration of $\mathcal{A}_k$. To do this, let us first consider an $R$-algebra $ C $. From Isomorphism (\ref{Bastheorem}) one gets a map
$$ \mathcal{A}(C) \rightarrow \mathcal{A}'(C \otimes_R R'), $$
which further gives 
$$ \mathcal{A}(C) \rightarrow \mathcal{A}'(C \otimes_R R') \rightarrow \mathcal{A}'(C \otimes_R R'/(\pi'^i)), $$
for any integer $i$ such that $ 0 \leq i \leq n $. Define functors $ F^i\mathcal{A}_k $ by
$$ F^i\mathcal{A}_k(C) = \Ker(\mathcal{A}(C) \rightarrow \mathcal{A}'(C \otimes_R R'/(\pi'^i))), $$
for any $ k = R/(\pi)$-algebra $C$. The functors $ F^i\mathcal{A}_k $ are represented by closed subgroup schemes of $ \mathcal{A}_k $, and give rise to a descending filtration

\begin{equation}\label{filtr}
\mathcal{A}_k = F^0\mathcal{A}_k \supseteq F^1\mathcal{A}_k \supseteq \ldots \supseteq F^n\mathcal{A}_k = 0. 
\end{equation}

\subsection{} 

The successive quotients of Filtration (\ref{filtr}) can be described quite accurately: Let $ Gr^i \mathcal{A}_k $ denote the quotient $ F^i\mathcal{A}_k/F^{i+1}\mathcal{A}_k $, for $ i \in \{ 0, \ldots, n-1 \} $. Then, according to Theorem 5.3 in \cite{Edix}, we have that $ Gr^0(\mathcal{A}_k) = (\mathcal{A}'_k)^{\boldsymbol{\mu}_n} $, and for $ 0 < i < n $, we have that
$$ Gr^i \mathcal{A}_k \cong T_{\mathcal{A}'_k,0}[i] \otimes_k (\mathfrak{m}/\mathfrak{m}^2)^{\otimes i}, $$
where $ \mathfrak{m} \subset R'$ is the maximal ideal, and where $ T_{\mathcal{A}'_k,0}[i] $ denotes the subspace of $ T_{\mathcal{A}'_k,0} $ where $ \xi \in \boldsymbol{\mu}_n $ acts by multiplication by $ \xi^i $. In particular, we note that the group schemes $ F^i\mathcal{A}_k $ are \emph{unipotent} for $ i > 0 $.

The filtration \emph{jumps} at the index $ i \in \{ 0, \ldots, n-1 \} $ if $ Gr^i \mathcal{A}_k \neq 0 $. Since 
$$ T_{\mathcal{A}'_k,0}[0] = (T_{\mathcal{A}'_k,0})^{\boldsymbol{\mu}_n} = T_{(\mathcal{A}'_k)^{\boldsymbol{\mu}_n},0} $$
(use \cite{Edix}, Proposition 3.2), it follows that the jumps are completely determined by the representation of $ \boldsymbol{\mu}_n $ on $ T_{\mathcal{A}'_k,0} $. In particular, it follows that there are at most $ \Dim(A) $ jumps, since $ \Dim_k T_{\mathcal{A}'_k,0} = \Dim(A) $.

\subsection{Filtration with rational indices}\label{ratfil}

Let $ a \in \mathbb{Z}_{(p)} \cap [0,1] $. If $ a = i/n $, then we define $ \mathcal{F}^a \mathcal{A}_k = F^i_n \mathcal{A}_k $, where $ F^i_n \mathcal{A}_k $ denotes the $i$-th step in the filtration induced by the tame extension of degree $n$. This definition does not depend on the choice of representatives $i$ and $n$ for $ a = i/n $ (\cite{Thesis}, Lemma 2.3). The following proposition is immediate: 

\begin{prop}\label{ratfilprop}
The construction above gives a descending filtration
$$ \mathcal{A}_k = \mathcal{F}^0 \mathcal{A}_k \supseteq \ldots \supseteq \mathcal{F}^a\mathcal{A}_k \supseteq \ldots \supseteq \mathcal{F}^1 \mathcal{A}_k = 0 $$
of $ \mathcal{A}_k $ by closed subgroup schemes, where $ a \in \mathbb{Z}_{(p)} \cap [0,1] $.
\end{prop}


\vspace{0.3cm}

 Let $ x \in [0,1] $ be a real number, and let $ (x^j)_j $ (resp. $ (x_k)_k $) be a sequence of numbers in $ \mathbb{Z}_{(p)} \cap [0,1] $ converging to $x$ from above (resp. from below). We will say that $ \{ \mathcal{F}^a \mathcal{A}_k \} $ \emph{jumps} at $x$ if $ \mathcal{F}^{x_k} \mathcal{A}_k \supsetneq \mathcal{F}^{x^j} \mathcal{A}_k $ for all $j$ and $k$. It is natural to ask \emph{how} many jumps there are, and \emph{where} they occur. It is easily seen that since every discrete     filtration $ \{ F^i_n \mathcal{A}_k \} $ jumps at most $ g = \Dim(A) $ times, the filtration $ \{ \mathcal{F}^a \mathcal{A}_k \} $ can have at most $g$ jumps. 

Consider a positive integer $n$ that is not divisible by $p$, and let $ \{ F^i_n \mathcal{A}_k \} $ be the filtration induced by the extension of degree $n$. Let us assume that this filtration has a jump at some $ i \in \{ 0, \ldots, n-1 \} $. Then we can say that $ \{ \mathcal{F}^a \mathcal{A}_k \} $ has a jump in the interval $ [i/n, (i+1)/n] $. By computing jumps in this way for increasing $n$, we get finer partitions of the interval $ [0,1] $, and increasingly better approximations of the jumps in $ \{ \mathcal{F}^a \mathcal{A}_k \} $. 

It follows that one can compute the jumps of $ \{ \mathcal{F}^a \mathcal{A}_k \} $ by computing the jumps for the filtrations $ \{ F^i_n \mathcal{A}_k \} $ for ``sufficiently'' many $n$ that are not divisible by $p$. This would for instance be the case for a multiplicatively closed subset $ \mathcal{U} \subset \mathbb{N} $ such that $ \mathbb{Z}[\mathcal{U}^{-1}] \cap [0,1] $ is dense in $ \mathbb{Z}_{(p)} \cap [0,1] $. 

\subsection{}
In the case where $ A/K $ obtains semi-abelian reduction over a tamely ramified extension $K'$ of $K$, the jumps of $ \{ \mathcal{F}^a \mathcal{A}_k \} $ have an interesting interpretation, which we will now explain. Let $ \widetilde{K} $ be the minimal extension over which $A$ aquires semi-abelian reduction (cf. ~\cite{Deschamps}, Th\'eor\`eme 5.15), and let $ \tilde{n} = \Deg(\widetilde{K}/K) $. Then we shall see below that the jumps occur at rational numbers of the form $ k/\tilde{n} $, where $ k \in \{0, \ldots, \tilde{n} - 1 \} $. This is essentially due to the following observation:

\begin{lemma}\label{tame jumps}
Let $ \widetilde{K}/K $ be the minimal extension over which $A/K$ obtains semi-abelian reduction, and let $ \tilde{n} = \Deg(\widetilde{K}/K) $. Consider a tame extension $ K'/K $ of degree $ n $, factoring via $ \widetilde{K} $, and let $ m = n/\tilde{n} $. Let $\mathcal{A}'/S' $ be the N\'eron model of $A_{K'}$. 

Then we have that the jumps in the filtration $ \{ F^i_n \mathcal{A}_k \} $ induced by $ S'/S $ occur at indices $ i = k n / \tilde{n} $, where $ 0 \leq k \leq \tilde{n} - 1$.
\end{lemma}
\begin{pf}
Let $ \widetilde{\mathcal{A}}/\widetilde{S} $ be the N\'eron model of $ A_{\widetilde{K}} $. By assumption, both $\mathcal{A}'$ and $ \widetilde{\mathcal{A}} $ are semi-abelian. Since $ \widetilde{\mathcal{A}}_{S'} $ is smooth, and $\mathcal{A}'$ has the N\'eronian property, we get a canonical morphism $ \widetilde{\mathcal{A}}_{S'} \rightarrow \mathcal{A}' $, extending the identity map on the generic fibers. Since $ \widetilde{\mathcal{A}}_{S'} $ is semi-abelian, it follows from Proposition 7.4/3 in \cite{Ner} that this morphism induces an isomorphism $ (\widetilde{\mathcal{A}}_k)^0 \cong (\mathcal{A}'_k)^0 $. In particular, we get that $ T_{\widetilde{\mathcal{A}}_k,0} = T_{\mathcal{A}'_k, 0} $.

Consider now the filtration $ \{ F^i_m \widetilde{\mathcal{A}}_k \} $ of $ \widetilde{\mathcal{A}}_k $ induced by the extension $ S'/\widetilde{S} $. Since $ \widetilde{\mathcal{A}} $ is semi-abelian, we have that $ F^i_m \widetilde{\mathcal{A}}_k = 0 $ for all $ i > 0 $. Therefore, we get that
$$ \widetilde{\mathcal{A}}_k = F^0_m \widetilde{\mathcal{A}}_k = Gr^0_m \widetilde{\mathcal{A}}_k = (\mathcal{A}'_k)^{\boldsymbol{\mu}_m}. $$
But now 
$$ (T_{\mathcal{A}'_k, 0})^{\boldsymbol{\mu}_m} = T_{(\mathcal{A}'_k)^{\boldsymbol{\mu}_m}, 0} = T_{\mathcal{A}'_k, 0}, $$ 
and so it follows that $ \boldsymbol{\mu}_m $ acts trivially on $ T_{\mathcal{A}'_k, 0} $. 

Let us now consider the filtration $ \{ F^i_n \mathcal{A}_k \} $ induced by the extension $ S'/S $. The jumps in this filtration are determined by the $ \boldsymbol{\mu}_n $-action on $ T_{\mathcal{A}'_k, 0} $. Assume that $ T_{\mathcal{A}'_k, 0}[i] \neq 0 $, for some $ i \in \{0, \ldots, n -1 \} $. On this subspace, every $ \xi \in \boldsymbol{\mu}_n $ acts by multiplication by $ \xi^i $. We can identify $ \boldsymbol{\mu}_m $ with the $\tilde{n}$-th powers in $ \boldsymbol{\mu}_n $, and since we established above that $ \boldsymbol{\mu}_m $ acts trivially, it follows that $ \xi^{\tilde{n} i} = 1 $. So therefore $ \tilde{n} i = k n $ for some $ k \in \{0, \ldots, \tilde{n} - 1 \} $, and we get that $ i = k n/\tilde{n} $.
\end{pf}

\vspace{0.5cm}

Using this lemma, one gets the following result:

\begin{prop}\label{tamejumpprop}
If $ A/K $ obtains semi-abelian reduction over a tamely ramified extension of $K$, then the jumps in the filtration $ \{ \mathcal{F}^a \mathcal{A}_k \} $ occur at indices $ k/\tilde{n} $, where $ k \in \{0, \ldots, \tilde{n} - 1 \} $, and where $ \tilde{n} $ is the degree of the minimal extension $ \widetilde{K}/K $ that realizes semi-abelian reduction for $ A $.
\end{prop}
\begin{pf}
Let us consider the sequence of integers $ (\tilde{n} m)_m $, where $m$ runs over the positive integers that are not divisibe by $p$. For the extension of degree $n = \tilde{n} m$, Lemma \ref{tame jumps} gives that the jumps of $ \{ F^i_n \mathcal{A}_k \} $ occur at indices $ i = k n / \tilde{n} $, where $ 0 \leq k \leq \tilde{n} - 1$. It follows that the jumps of $ \{ \mathcal{F}^a \mathcal{A}_k \} $ will be among the limits of the expressions $ i/n = k/ \tilde{n} $, as $m$ goes to infinity, and the result follows.
\end{pf}

\vspace{0.3cm}

The next proposition shows that the jumps come in ''simultaneously reduced'' form.

\begin{prop}\label{minimalityjump}
Let us assume that $ A/K $ obtains semi-abelian reduction over a tamely ramified minimal extension $ \widetilde{K}/K$, and that $ \tilde{n} = [\widetilde{K}:K] > 1 $. Let $ i_1/\tilde{n}, \ldots, i_g/\tilde{n} $ be the jumps in the filtration $ \{ \mathcal{F}^a \mathcal{A}_k \} $.

Assume that $m$ is a positive integer such that $ m $ divides $ i_l $ for all $l$, and that $ m $ divides $ \tilde{n} $. Then it follows that $ m = 1 $.
\end{prop}
\begin{pf}
Let us assume to the contrary that $ m > 1 $. Then it follows that the subgroup $ H \subseteq \boldsymbol{\mu}_{\tilde{n}} $ consisting of $\tilde{n}/m $-th powers acts trivially on $ T_{\widetilde{\mathcal{A}}_k, 0} $. We now claim that $H$ acts trivially also on $ \widetilde{\mathcal{A}}_k^0 $. To see this, we first observe that 
$$ T_{\widetilde{\mathcal{A}}_k, 0} \cong (T_{\widetilde{\mathcal{A}}_k, 0})^H \cong T_{\widetilde{\mathcal{A}}_k^H, 0}. $$
Furthermore, since $p$ does not divide the order of $H$, we have that $\widetilde{\mathcal{A}}_k^H$ is smooth, and that the canonical inclusions $ \widetilde{\mathcal{A}}_k^H \subseteq \widetilde{\mathcal{A}}_k $ and $ (\widetilde{\mathcal{A}}_k^H)^0 \subseteq \widetilde{\mathcal{A}}_k^0 $ are closed immersions. Hence, $\Dim((\widetilde{\mathcal{A}}_k^H)^0) = \Dim((\widetilde{\mathcal{A}}_k)^0) $, and so it follows that $ (\widetilde{\mathcal{A}}_k^H)^0 = (\widetilde{\mathcal{A}}_k^0)^H = \widetilde{\mathcal{A}}_k^0 $. Therefore, $H$ acts trivially on $\widetilde{\mathcal{A}}_k^0$. But, by Lemma 5.16 in \cite{Deschamps}, this contradicts the minimality of $ \widetilde{K}/K $. 
\end{pf}

\subsection{The case of Jacobians}\label{jacobiancase}
Let $X/K$ be a smooth, projective and geometrically connected curve of genus $ g > 0 $. Let $ J' = J_{K'} $ denote the Jacobian of $X_{K'}$, and let $ \mathcal{J}'/S' $ be the N\'eron model of $J'$ over $S'$. 

We can let $G$ act on $ X_{K'} $ via the action on the second factor. Let $ \mathcal{Y}/S' $ be a regular model of $ X_{K'} $ such that the $G$-action on $ X_{K'} $ extends to $ \mathcal{Y} $. According to \cite{Ner}, Theorem 9.5/4, there is, under certain hypotheses, a canonical isomorphism
$$ \Pic^0_{\mathcal{Y}/S'} \cong \mathcal{J}'^0, $$
where $ \mathcal{J}'^0 $ is the identity component of $ \mathcal{J}' $, and where $ \Pic^0_{\mathcal{Y}/S'} $ is the identity component of the relative Picard functor $ \Pic_{\mathcal{Y}/S'} $. Hence, on the special fibers, we get an isomorphism
$$ \Pic^0_{\mathcal{Y}_k/k} \cong \mathcal{J}'^0_k. $$
By \cite{Ner}, Theorem 8.4/1, it follows that we can canonically identify
\begin{equation}\label{H^1=T}
H^1(\mathcal{Y}_k, \mathcal{O}_{\mathcal{Y}_k}) \cong T_{\mathcal{J}'_k,0}. 
\end{equation}

We are interested in computing the irreducible characters for the representation of $\boldsymbol{\mu}_n$ on $T_{\mathcal{J}'_k,0}$. With the identification in \ref{H^1=T} above, we see that this can be done by computing the irreducible characters for the representation of $\boldsymbol{\mu}_n$ on $ H^1(\mathcal{Y}_k, \mathcal{O}_{\mathcal{Y}_k}) $.

By combining the discussion in this section with properties of the representation of $\boldsymbol{\mu}_n$ on $ H^1(\mathcal{Y}_k, \mathcal{O}_{\mathcal{Y}_k}) $, we obtain in Corollary \ref{specific jump corollary} a quite precise description of the jumps of the filtration $ \{ \mathcal{F}^a \mathcal{J}_k \} $.

\section{Tame extensions and Galois actions}\label{extensions and actions}

\subsection{}
Throughout the rest of this paper, $R$ will denote a complete discrete valuation ring, with fraction field $K$, and with algebraically closed residue field $k$. 

$X/K$ will be a smooth, projective, geometrically connected curve over $K$, of genus $g(X)>0$.

\begin{dfn}
A scheme $ \mathcal{X} $ is called a \emph{model} of $X$ over $ S = \Spec(R) $ if $ \mathcal{X} $ is integral and normal, projective and flat over $S$, and with generic fibre $ \mathcal{X}_K \cong X $.
\end{dfn}



It is well known that we can always find a regular model for $X/K$ (see for instance \cite{Liubook}). By blowing up points in the special fiber, we can even ensure that the irreducible components of the special fiber are smooth, and intersect transversally. Such a model will be called a \emph{strict normal crossings} model for $X/K$, or for short, an SNC-model.

\subsection{}\label{2.2}
Let $ \mathcal{X}/S $ be an SNC-model for $X/K$. Let $ K \subset K' $ be a finite, separable field extension, and let $R'$ be the integral closure of $R$ in $K'$. Since $R$ is complete, we have that $R'$ is a complete discrete valuation ring (\cite{Serre}, Proposition II.3). Making the finite base extension $S' = \Spec(R') \rightarrow S = \Spec(R) $,  we obtain a commutative diagram
$$ \xymatrix{
 \mathcal{Y} \ar[d] \ar[r]^{\rho} & \mathcal{X}' \ar[d] \ar[r]^{f}  & \mathcal{X} \ar[d]\\ 
  S' \ar[r]^{\Id} & S' \ar[r] & S, } $$
where $ \mathcal{X}' $ is the normalization of the pullback $ \mathcal{X}_{S'} = \mathcal{X} \times_S S' $ ($ \mathcal{X}_{S'} $ is integral by Lemma \ref{lemma construction} below), and $ \rho : \mathcal{Y} \rightarrow \mathcal{X}' $ is the minimal desingularization. The map $ f : \mathcal{X}' \rightarrow \mathcal{X} $ is the composition of the projection $ \mathcal{X}_{S'} \rightarrow \mathcal{X} $ with the normalization $ \mathcal{X}' \rightarrow \mathcal{X}_{S'} $.

\begin{lemma}\label{lemma construction} With the hypotheses above, the following statements hold:
\begin{enumerate}
\item The pullback $ \mathcal{X}_{S'} $ is an integral scheme.
\item $ f : \mathcal{X}' \rightarrow \mathcal{X} $ is a finite morphism.
\end{enumerate}
\end{lemma}
\begin{pf}
(i) Let us first note that the generic fiber of $ \mathcal{X}_{S'} $ is the pullback $\mathcal{X}_K \otimes_K K' $, where $ \mathcal{X}_K $ is the generic fiber of $ \mathcal{X} $. By assumption $ \mathcal{X}_K $ is smooth and geometrically connected over $K$, so in particular the generic fiber of $ \mathcal{X}_{S'} $ is integral. Now, since $ \mathcal{X}_{S'} \rightarrow S' $ is flat, it follows from \cite{Liubook}, Proposition 4.3.8, that $ \mathcal{X}_{S'} $ is integral as well.

(ii) Since $ R' $ is a \emph{complete} discrete valuation ring it is excellent (\cite{Liubook}, Theorem 8.2.39). As $ \mathcal{X}_{S'} $ is of finite type over $S'$, it follows that $ \mathcal{X}_{S'} $ is an excellent scheme, and hence the normalization morphism $ \mathcal{X}' \rightarrow \mathcal{X}_{S'} $ is finite (\cite{Liubook}, Theorem 8.2.39). The projection $ \mathcal{X}_{S'} \rightarrow \mathcal{X} $ is finite, since it is the pullback of the finite morphism $ S' \rightarrow S $. So the composition $f$ of these two morphisms is indeed finite.
\end{pf}

\subsection{} Let us now assume that the field extension $ K \subset K' $ is Galois with group $G$. Every $ \sigma \in G $ induces an automorphism of $R'$ that fixes $R$, and we have furthermore that $ R'^{G} = R $. So there is an injective group homomorphism $ G \rightarrow \Aut(S') $, and we may view $ S' \rightarrow S $ as the quotient map. 

We can lift the $G$-action to $ \mathcal{X}_{S'} $, via the action on the second factor. So there is a group homomorphism $ G \rightarrow \Aut(\mathcal{X}_{S'}) $. For any element $ \sigma \in G $, we shall still denote the image in $ \Aut(\mathcal{X}_{S'}) $ by $ \sigma $. Proposition \ref{prop. 2.3} below states that this action lifts uniquely both to the normalization $ \mathcal{X}' $ and to the minimal desingularization $ \mathcal{Y} $ of $ \mathcal{X}' $.


\begin{prop}\label{prop. 2.3} With the hypotheses above, the following statements hold:
\begin{enumerate}
\item The $G$-action on $ \mathcal{X}_{S'} $ lifts uniquely to the normalization $ \mathcal{X}' $.

\item The $ G $-action on $ \mathcal{X}' $ lifts uniquely to the minimal desingularization $ \mathcal{Y} $.

\item For any $ \sigma \in G $, let $ \sigma $ denote the induced automorphism of $ \mathcal{X}' $, and let $ \tau $ be the unique lift of $ \sigma $ to $ \Aut(\mathcal{Y}) $. Then we have that $ \tau(\rho^{-1}(\Sing(\mathcal{X}'))) = \rho^{-1}(\Sing(\mathcal{X}')) $. That is, the exceptional locus is mapped into itself under the $G$-action on $ \mathcal{Y} $.
\end{enumerate}
\end{prop}
\begin{pf}
This is straightforward from the universal properties of the normalization and of the minimal desingularization. For a detailed proof, we refer to \cite{Thesis}.
\end{pf}

\subsection{}\label{assumption on degree}
We shall throughout the rest of the paper make the assumption that $ n = [K':K] $ is not divisible by the residue characteristic $p$. Since $k$ is algebraically closed it has a full set $ \boldsymbol{\mu}_n $ of $n$-th roots of unity, and as $R$ is complete, we may lift all $n$-th roots of unity to $R$. We can choose a uniformizing parameter $ \pi \in R $ such that $ K' = K[\pi']/(\pi'^n - \pi) $. The extension $ K \subset K' $ is Galois, with group $ G = \boldsymbol{\mu}_n $. Also, $ R' := R[\pi']/(\pi'^n - \pi) $ is the integral closure of $R$ in $K'$, and $ \pi'$ is a uniformizing parameter for $R'$.

\subsection{Assumptions on $\mathcal{X}$}\label{assumption on surface}
Throughout the rest of this paper, we shall make two assumptions in the situation considered in Section \ref{2.2}:

\begin{ass}\label{ass. 2.4}
Let $ x \in \mathcal{X} $ be a closed point in the special fiber such that two irreducible components $C_1$ and $C_2$ of $ \mathcal{X}_k $ meet at $x$, and let $ m_i = \Mult(C_i) $. We will always assume that \emph{at least} one of the $m_i$ is not divisible by $p$. 
\end{ass}

With this assumption, we can find an isomorphism 
$$ \widehat{\mathcal{O}}_{\mathcal{X},x} \cong R[[u_1,u_2]]/(\pi - u_1^{m_1} u_2^{m_2}) $$
(cf. \cite{CED}, proof of Lemma 2.3.2).

\begin{ass}\label{ass. 2.5}
Let $l$ denote the least common multiple of the multiplicities of the irreducible components of $\mathcal{X}_k$. Then we assume that $ \Gcd(l,n) = 1 $. 
\end{ass}

When Assumptions \ref{ass. 2.4} and \ref{ass. 2.5} are valid, the following facts can be proved using the computations in \cite{Thesis}:

$\bullet$ Let $ x \in \mathcal{X} $ be a closed point in the special fiber. Because of Assumption \ref{ass. 2.5}, there is a unique point $ x' \in \mathcal{X}'_k $ that maps to $x$. The local analytic structure of $ \mathcal{X}' $ at $x'$ depends only on $n = [K':K]$ and on the local analytic structure of $ \mathcal{X} $ at $x$. If $x$ belongs to a unique component of $ \mathcal{X}_k $, then $x'$ belongs to a unique component of $ \mathcal{X}'_k $, and $ \mathcal{X}' $ is regular at $x'$. If $x$ is an intersection point of two distinct components, then the same is true for $x'$, and $ \mathcal{X}' $ will have a \emph{tame cyclic quotient singularity} at $x'$. 

$\bullet$ The minimal desingularization $ \mathcal{Y} $ of $ \mathcal{X}' $ is an SNC-model. Furthermore, the structure of $ \mathcal{Y} $ locally above a tame cyclic quotient singularity $x' \in \mathcal{X}'$ is completely determined by the structure locally at $ x = f(x') $ and the degree $n$ of the extension. The inverse image of $x'$ consists of a chain of smooth and rational curves whose multiplicities and self intersection numbers may be computed from the integers $n, m_1$ and $m_2$.

$\bullet$ For every irreducible component $C$ of $ \mathcal{X}_k $, there is precisely one component $C'$ of $ \mathcal{X}'_k $ that dominates $C$. The component $C'$ is isomorphic to $C$, and we have that $ \Mult_{\mathcal{X}'_k}(C') = \Mult_{\mathcal{X}_k}(C) $. It follows that the combinatorial structure of $ \mathcal{Y}_k $ is completely determined by the combinatorial structure of $ \mathcal{X}_k $ and the degree of $S'/S$.

\subsection{}
We will now begin to describe the $G$-action on $ \mathcal{X}' $ and $ \mathcal{Y} $ in more detail. Assumptions \ref{ass. 2.4} and \ref{ass. 2.5} will impose some restrictions on this action.

\begin{prop}\label{action on desing}
Let $ \rho : \mathcal{Y} \rightarrow \mathcal{X'} $ be the minimal desingularization. Then the following properties hold: 
\begin{enumerate}
\item Let $D$ be an irreducible component of $ \mathcal{Y}_k $ that dominates a component of $ \mathcal{X}_k $. Then $D$ is stable under the $G$-action, and $G$ acts \emph{trivially} on $D$.
\item Let $ x' \in \mathcal{X}' $ be a singular point, and let $E_1, \ldots, E_{l}$ be the exceptional components mapping to $x'$ under $ \rho $. Then every $E_i$ is stable under the $G$-action, and every node in the chain $ \rho^{-1}(x') $ is fixed under the $G$-action.
\end{enumerate}
\end{prop}
\begin{pf}
Let us first note that the map $ \mathcal{X}_{S'} \rightarrow \mathcal{X} $ is an isomorphism on the special fibers. Moreover, the action on the special fiber of $ \mathcal{X}_{S'} $ is easily seen to be trivial, so every closed point in the special fiber is fixed. Since the action on $ \mathcal{X}' $ commutes with the action on $ \mathcal{X}_{S'} $, it follows that every point in the special fiber of $ \mathcal{X}' $ is fixed. In particular, every irreducible component $C'$ of $ \mathcal{X}'_k $ is stable under the $G$-action, and the restriction of this action to $C'$ is trivial. Since the action on $ \mathcal{Y} $ commutes with the action on $ \mathcal{X}' $, it follows that the same is true for the strict transform $D$ of $C'$ in $ \mathcal{Y} $. This proves (i).

For (ii), we observe that since $x'$ is fixed, we have that $\rho^{-1}(x')$ is stable under the $G$-action. But also the two branches meeting at $x'$ are fixed. Let $D$ be the strict transform of any of these two branches. From part (i), it follows that the point where it meets the exceptional chain $\rho^{-1}(x')$ must be fixed. So if $E_1$ is the component in the chain meeting $\widetilde{D}$, then $E_1$ must be mapped into itself. Let $E_2$ be the next component in the chain. Then the point where $E_1$ and $E_2$ meet must also be fixed, so $E_2$ must also be mapped to itself. Continuing in this way, it is easy to see that all of the exceptional components are stable under the $G$-action, and that all nodes in $\rho^{-1}(x')$ are fixed points.  
\end{pf}

\begin{cor}\label{g^{-1}(Z) = Z}
Let $ 0 \leq Z \leq \mathcal{Y}_k $ be an effective divisor. Then the $G$-action restricts to $Z$.
\end{cor}
\begin{pf}
 Since $Z$ is an effective Weil divisor, we can write $ Z = \sum_{C} r_C C $, where $C$ runs over the irreducible components of $ \mathcal{Y}_k $, and $r_C$ is a non-negative integer for all $C$. But Proposition \ref{action on desing} states that all irreducible components $C$ of $ \mathcal{Y}_k $ are stable under the $G$-action, and hence we get that the same holds for $ Z $. In other words, the action restricts to $Z$.
\end{pf}

\vspace{0.3cm}

From Proposition \ref{action on desing} above, it follows that every node $y$ in $ \mathcal{Y}_k $ is a fixed point for the $G$-action on $ \mathcal{Y} $. Hence there is an induced action on $ \mathcal{O}_{\mathcal{Y},y} $ and on the cotangent space $ \mathfrak{m}_y/\mathfrak{m}_y^2 $, where $ \mathfrak{m}_y \subset \mathcal{O}_{\mathcal{Y},y} $ is the maximal ideal. In order to get a precise description of the action on the cotangent space, we will first describe the action on the completion $ \widehat{\mathcal{O}}_{\mathcal{Y},y} $.
 
Since, by Proposition \ref{action on desing}, every irreducible component $D$ of $ \mathcal{Y}_k $ is mapped to itself under the $G$-action, it follows that the $G$-action restricts to $D$ and that the points where $D$ meets the rest of the special fiber are fixed. In the case where $G$ acts non-trivially on $D$, we will see in Proposition \ref{prop. 3.3} that the fixed points for the $G$-action on $D$ are precisely the points where $D$ meets the rest of the special fiber. In particular, we wish to describe the action on $D$ locally at the fixed points. 

\section{Desingularizations and actions}\label{lifting the action}

In this section, we study how one can explicitly describe the action on the minimal desingularization $ \rho : \mathcal{Y} \rightarrow \mathcal{X}' $. Since we are only interested in this action locally at fixed points or stable components in the exceptional locus of $\rho$, we will begin with showing that we can reduce to studying the minimal desingularization locally at a singular point $ x' \in \mathcal{X}' $. This is an important step, since we have a good description of the complete local ring $ \widehat{\mathcal{O}}_{\mathcal{X}',x'} $. In particular, we can find a nice algebraization of this ring, with a compatible $G$-action. It turns out that it suffices for our purposes to study the minimal desingularization of this ring, and the lifted $G$-action.

In the second part of this section, we study the desingularization of an algebraization of $ \widehat{\mathcal{O}}_{\mathcal{X}',x'} $. We use the explicit blow up procedure in \cite{CED} in order to describe how the $G$-action lifts. In particular, we describe the action on the completion of the local rings at the nodes in the exceptional locus, and the action on the exceptional components. These results are collected in Proposition \ref{prop. 3.3}.   

\subsection{}\label{seclocred}
If $ x' \in \mathcal{X}' $ is a singular point, we need to understand how $G$ acts on $ \widehat{\mathcal{O}}_{\mathcal{X}',x'} $. In order to do this, we consider the image $ f(x') = x $ of $x'$ under the morphism $ f : \mathcal{X}' \rightarrow \mathcal{X} $. Then $x$ is a closed point in the special fiber, and we have that   
$$ \widehat{\mathcal{O}}_{\mathcal{X},x} \cong R[[v_1,v_2]]/(\pi - v_1^{m_1} v_2^{m_2}), $$
where $m_1$ and $m_2$ are positive integers. Let $n$ be the degree of $ R'/R $, which by assumption is relatively prime to $m_1$ and $m_2$. In the discussion that follows we will use some properties that were proved in \cite{Thesis}.

We let $ G = \boldsymbol{\mu}_n $ act on $ \mathcal{X}_{S'} $ via its action on the second factor. We point out that we here choose the action on $R'$ given by $ [\xi](\pi') = \xi \pi' $ for any $ \xi \in \boldsymbol{\mu}_n $. Choosing this action is notationally convenient when we work with rings. However, the natural right action on $ \mathcal{X}_{S'} $ is the inverse to the one we use here. In particular, the irreducible characters for the representation of $ \boldsymbol{\mu}_n $ on $ H^1(\mathcal{Y}_k, \mathcal{O}_{\mathcal{Y}_k}) $ induced by the action chosen here on $\mathcal{X}_{S'}$ will be the inverse characters to those induced by the right action.

Let also $ x $ denote the unique point of $ \in \mathcal{X}_{S'} $ mapping to $ x \in \mathcal{X} $. The map $ \mathcal{O}_{\mathcal{X},x} \rightarrow \mathcal{O}_{\mathcal{X}_{S'},x} $ associated to the projection $ \mathcal{X}_{S'} \rightarrow \mathcal{X} $ can be described by the tensorization
$$ \mathcal{O}_{\mathcal{X},x} \rightarrow \mathcal{O}_{\mathcal{X},x} \otimes_R R', $$
and the $G$-action on $ \mathcal{O}_{\mathcal{X}_{S'},x} = \mathcal{O}_{\mathcal{X},x} \otimes_R R' $ is induced from the action on $R'$.

Since $ \mathcal{O}_{\mathcal{X},x} \rightarrow \mathcal{O}_{\mathcal{X},x} \otimes_R R' $ is finite, completion commutes with tensoring with $R'$, so we get that
$$ \widehat{\mathcal{O}}_{\mathcal{X}_{S'},x} = \widehat{\mathcal{O}}_{\mathcal{X},x} \otimes_R R', $$
and hence the $G$-action on $ \widehat{\mathcal{O}}_{\mathcal{X}_{S'},x} $ is induced from the action on $R'$ in the second factor. It follows that
$$ \widehat{\mathcal{O}}_{\mathcal{X}_{S'},x} \cong R'[[v_1,v_2]]/(\pi'^n - v_1^{m_1} v_2^{m_2}), $$
and that the $G$-action can be described by $ [\xi](\pi') = \xi \pi' $ and $ [\xi](v_i) = v_i $, for any $ \xi \in \boldsymbol{\mu}_n $.

Let $ \mathcal{X}' \rightarrow \mathcal{X}_{S'} $ be the normalization. There is a unique point $x'$ mapping to $x$, and the induced map $ \widehat{\mathcal{O}}_{\mathcal{X}_{S'},x} \rightarrow \widehat{\mathcal{O}}_{\mathcal{X}',x'} $ is the normalization of $ \widehat{\mathcal{O}}_{\mathcal{X}_{S'},x} $. Furthermore, the $G$-action on $ \widehat{\mathcal{O}}_{\mathcal{X}',x'} $ induced by the action on $ \mathcal{X}' $ is the unique lifting of the $G$-action on $ \widehat{\mathcal{O}}_{\mathcal{X}_{S'},x} $ to the normalization $ \widehat{\mathcal{O}}_{\mathcal{X}',x'} $.

Let $ \rho : \mathcal{Y} \rightarrow \mathcal{X}' $ be the minimal desingularization, and consider the fiber diagram
$$ \xymatrix{
\widehat{\mathcal{Y}} \ar[d]_{\hat{\rho}} \ar[r]^{\phi} & \mathcal{Y} \ar[d]^{\rho} \\
\Spec(\widehat{\mathcal{O}}_{\mathcal{X}',x'}) \ar[r] & \mathcal{X}' .}
$$
Then $ \hat{\rho} $ is the minimal desingularization of $ \Spec(\widehat{\mathcal{O}}_{\mathcal{X}',x'}) $ (cf. ~\cite{Lip}, Lemma 16.1, and use the fact that $ \mathcal{Y} $ is minimal), and hence the $G$-action on $ \Spec(\widehat{\mathcal{O}}_{\mathcal{X}',x'}) $ lifts uniquely to $ \widehat{\mathcal{Y}} $. 

The projection $\phi$ induces an isomorphism of the exceptional loci $ \hat{\rho}^{-1}(x') $ and $ \rho^{-1}(x') $. Let $E$ be an exceptional component. Then the $G$-action restricts to $E$, and it is easily seen that $\phi$, when restricted to $E$, is equivariant. 

Furthermore, for any closed point $ y \in \rho^{-1}(x') $, we have that $\phi$ induces an isomorphism $ \widehat{\mathcal{O}}_{\mathcal{Y},y} \cong \widehat{\mathcal{O}}_{\widehat{\mathcal{Y}},y} $ (one can argue in a similar way as in the proof of \cite{Liubook}, Lemma 8.3.49). If $y$ is a fixed point, it is easily seen that this isomorphism is equivariant.
We therefore conclude that in order to describe the action on $\mathcal{Y}$ locally at the exceptional locus over $x'$, it suffices to consider the minimal desingularization of $ \Spec(\widehat{\mathcal{O}}_{\mathcal{X}',x'}) $.

\subsection{}

In order to find an algebraization of $ \widehat{\mathcal{O}}_{\mathcal{X}',x'} $, we consider first the polynomial ring $ V = R'[v_1,v_2]/(\pi'^n - v_1^{m_1} v_2^{m_2}) $. We let $G$ act on $ V $ by $ [ \xi ](\pi') = \xi \pi' $ and $ [ \xi ](v_i) = v_i $ for $ i = 1,2 $, for any $ \xi \in G $. Note that the maximal ideal $ \mathfrak{p} = (\pi', v_1,v_2) $ is fixed, and hence there is an induced action on the completion $ \widehat{V}_{\mathfrak{p}} = R'[[v_1,v_2]]/(\pi'^n - v_1^{m_1} v_2^{m_2}) $, given as above. This gives a $G$-equivariant algebraization of $ \widehat{\mathcal{O}}_{\mathcal{X}_{S'},x} $. 

Consider the $R'$-algebra homomorphism
$$ V = R'[v_1,v_2]/(\pi'^n - v_1^{m_1} v_2^{m_2}) \rightarrow T = R'[t_1,t_2]/(\pi' - t_1^{m_1} t_2^{m_2}), $$
given by $ v_i \mapsto t_i^n $. We let $ \boldsymbol{\mu}_n $ act on $T$, relatively to $R'$, by $ [ \eta ](t_1) = \eta t_1 $, $ [ \eta ](t_2) = \eta^r t_2 $, where $ r $ is the unique integer $ 0 < r < n $ such that $ m_1 + r m_2 \equiv_n 0 $. (Note that this is an ad hoc action introduced to compute the normalization, which must not be confused with the natural $G$-action). Arguing as in \cite{Thesis}, one can show that the induced map $ V \rightarrow U := T^{\boldsymbol{\mu}_n} $ is the normalization of $V$. Furthermore, it is easily seen that there is a unique maximal ideal $\mathfrak{q}  \subset U $ mapping to $ \mathfrak{p} $, corresponding to the "origin" $ (\pi', t_1, t_2) $ in $T$.

It is shown in \cite{Thesis} Lemma 4.1 that $U$ is an equivariant algebraization of $ \widehat{\mathcal{O}}_{\mathcal{X}',x'} $. Let $ \rho_U : \mathcal{Z} \rightarrow \Spec(U) $ be the minimal desingularization. Then we have a fiber diagram
$$ 
\xymatrix{
\widehat{\mathcal{Y}} \ar[d]_{\hat{\rho}} \ar[r] & \mathcal{Z} \ar[d]^{\rho_U} \\ 
\Spec(\widehat{\mathcal{O}}_{\mathcal{X}',x'}) \ar[r] & \Spec(U),}
$$
where all maps commute with the various $G$-actions. We conclude, by similar argumentation as in Section \ref{seclocred}, that in order to describe the $G$-action on $ \widehat{\mathcal{Y}} $ locally at fixed points or components in the exceptional locus, it suffices to compute the corresponding data for $ \mathcal{Z} $.

\subsection{}\label{3.11}

Proposition \ref{prop. 3.3} below gives a description of the $G$-action on $ \mathcal{Z} $. Having this description will be important in later sections, when we consider the $G$-action on the cohomology groups $ H^i(\mathcal{Y}_k, \mathcal{O}_{\mathcal{Y}_k}) $. For a proof, we refer to \cite{Thesis}, Proposition 4.3. In order to state this result, and for future reference, we list some properties associated to the resolution of the singularity $ Z = \Spec(U) $ (see \cite{CED} for proofs).

We call the integers $m_1$, $m_2$ and $n$ the \emph{parameters} of the singularity. Let $r$ be the unique integer with $ 0<r<n $ such that $ m_1 + rm_2 = 0 $ modulo $n$. Write $ \frac{n}{r} = [b_1, \ldots , b_l, \ldots , b_L]_{JH} $ for the Jung-Hirzebruch expansion. The exceptional locus of $ \rho_U $ consists of a \emph{string} of smooth and rational curves $ C_1, \ldots, C_L $, with self intersection numbers $ C_l^2 = - b_l $ and multiplicities $\mu_l$, for all $ l \in \{1, \ldots, L \} $.  

There are two series of numerical equations associated to the singularity. We have
\begin{equation}\label{equation 8.1} 
r_{l-1} = b_{l+1} r_l - r_{l+1}, 
\end{equation}
for $ 0 \leq l \leq L-1 $, where we put $ r_{-1} = n $ and $ r_0 = r $. Furthermore, we have
 \begin{equation}\label{equation 8.2} 
\mu_{l+1} = b_l \mu_l - \mu_{l-1}, 
\end{equation}
which is valid for $ 1 \leq l \leq L $. Here we define $ \mu_0 = m_2 $ and $ \mu_{L+1} = m_1 $. 

We also have the equation $ m_1 + r m_2 = n \mu_1 $ (see \cite{CED}, Corollary 2.4.3). Together with Equations (\ref{equation 8.1}) and (\ref{equation 8.2}), this equation enables you to compute the branch multiplicities.


\begin{prop}\label{prop. 3.3}
The minimal desingularization $ \mathcal{Z} $ of $ Z = \Spec(U) $ can be covered by the affine charts $ \Spec(U_l) $, where 
$$ U_l = R'[z_{l-1},w_{l-1}]/(z_{l-1}^{\mu_{l}} w_{l-1}^{\mu_{l-1}} - \pi'), $$
for $ l \in \{ 1, \ldots, L +1 \} $. 

These charts are $ G $-stable, and the $ G $-action is given by
$ [ \xi ] (\pi') = \xi \pi' $, $ [ \xi ] (z_{l-1}) = \xi^{ \alpha_1  r_{l-2} } z_{l-1} $ and $ [ \xi ] (w_{l-1}) = \xi^{ - \alpha_1 r_{l-1} } w_{l-1} $, where $\alpha_i$ denotes an inverse to $m_1$ modulo $n$, for all $ l \in \{ 1, \ldots, L+1 \} $, and for any $ \xi \in G $.

Let $ C_l $ be the $ l $-th exceptional component. On the chart $\Spec(U_l)$, we have that the affine ring for $ C_l $ is $ k[w_{l-1}] $, and $ G $ acts by $ [ \xi ] (w_{l-1}) = \xi^{ - \alpha_1 r_{l-1} } w_{l-1} $, for any $ \xi \in G $. On the chart $\Spec(U_{l+1})$, the affine ring for $ C_l $ is $ k[z_l] $, and $ G $ acts by $ [ \xi ] (z_{l}) = \xi^{ \alpha_1  r_{l-1} } z_{l} $.
\end{prop}

\vspace{0.3cm}

The following corollary is immediate from Proposition \ref{prop. 3.3}: 

\begin{cor}
The irreducible components $C_l$ of the exceptional locus are stable under the $G$-action. Furthermore, if $ \xi \in \boldsymbol{\mu}_n $ is a primitive root, then the automorphism of $C_l$ induced by $\xi$ is non-trivial, for all $ l \in \{ 1, \ldots, L \} $, with fixed points precisely at the two points where $C_l$ meets the rest of the special fiber.
\end{cor}

Let us finally remark that the cotangent space to $ \mathcal{Z} $ at the fixed point that is the intersection point of $ C_l $ and $ C_{l+1} $ is generated by (the classes) of the local equations $ z_l $ and $ w_l $ for the curves. Therefore, Proposition \ref{prop. 3.3} gives a complete description of the action on the cotangent space. Furthermore, we can also read off the eigenvalues for the elements of this basis. Hence we immediately get an explicit description of the action on the cotangent space to the minimal desingularization of $\mathcal{X}'$ at the corresponding fixed point. 

\begin{ex}\label{Gdesingex}
Consider the singularity with parameters $ (m_1,m_2,n) = (1,3,7) $. From the equation $ m_1 + r_0 m_2 = n \mu_1 $ we easily compute that $ r_0 = 2 $ and $ \mu_1 = 1 $. From the equation $ n = r_{-1} = b_1 r_0 - r_1 $ we find that $ b_1 = 4 $ and that $ r_1 = 1 $. So we get that $ L = 2 $, and hence there are two exceptional curves $C_1$ and $C_2$. In order to compute $\mu_2$, we use the equation $ \mu_0 + \mu_2 = b_1 \mu_1 $, and find that $ \mu_2 = 1 $.

We conclude this example with writing out the $G$-action on $C_1$, and the $G$-action on the cotangent space to $ \mathcal{Z} $ at the point $ y_1 = C_1 \cap C_2 $. In the notation of Proposition \ref{prop. 3.3}, $C_1$ and $C_2$ are generically contained in the $G$-stable open affine $ \Spec(U_2) $. We can immediately read off that the cotangent space is generated by (the classes of) $ z_1 $ and $ w_1 $. Furthermore, we have that the $ G $-action is given by $ [\xi](z_1) = \xi^{\alpha_1 r_0} z_1 = \xi^{2} z_1 $ and $ [\xi](w_1) = \xi^{- \alpha_1 r_1} w_1 = \xi^{6} w_1 $, for any $ \xi \in \boldsymbol{\mu}_{7} $.
 The affine ring for $C_1$ on this chart is $k[z_1]$, and the action is given by $ [\xi](z_1) = \xi^{2} z_1 $.
\end{ex}

\section{Computing traces}

We will now study of the $G$-action on the cohomology groups $ H^i(\mathcal{Y}_k, \mathcal{O}_{\mathcal{Y}_k}) $. Since the residue characteristic is (possibly) positive, we have to introduce the concept of \emph{Brauer characters}. This roughly amounts to lifting all eigenvalues to characteristic zero. 

A second problem is the fact that $ \mathcal{Y}_k $ in general is singular. Therefore, we shall first consider the case of a group acting on the cohomology groups of an invertible sheaf on a smooth projective curve. For such situations, one can write down trace formulas in terms of local data at the fixed points. Later, in Section \ref{section 6}, we shall use such computations in order to compute the characters for the action by $G$ on $ H^i(\mathcal{Y}_k, \mathcal{O}_{\mathcal{Y}_k}) $.




\subsection{}\label{inverse}
We begin this section with recalling some generalities about coherent sheaf cohomology, and establish notation and terminology that will be used throughout the rest of the paper. 

Let $ g : X \rightarrow Y $ be a morphism of schemes, and let $ \mathcal{F} $ be an $ \mathcal{O}_Y $-module. The morphism $ g $ induces a natural and canonical homomorphism
$$ H^p(g) : H^p(Y,\mathcal{F}) \rightarrow H^p(X, g^* \mathcal{F}), $$  
for all $p \geq 0$. 


Consider now the case where $Y = X$, so that $ g : X \rightarrow X $ is an endomorphism, and where $ \mathcal{F} $ is a sheaf of $ \mathcal{O}_X $-modules. Assume in addition that we are given a homomorphism $ u : g^* \mathcal{F} \rightarrow \mathcal{F} $ of $ \mathcal{O}_X $-modules. By functoriality, $u$ induces a homomorphism
$$ H^p(u) : H^p(X, g^* \mathcal{F}) \rightarrow H^p(X,\mathcal{F}), $$ 
for all $p \geq 0$.


\begin{dfn}\label{dfn 5.4}
Let $ g : X \rightarrow X $ be morphism, let $ \mathcal{F} $ be a sheaf of $ \mathcal{O}_X $-modules, and let $ u : g^* \mathcal{F} \rightarrow \mathcal{F} $ be a homomorphism of $ \mathcal{O}_X $-modules. The endomorphism 
$$ H^p(g,u) : H^p(X,\mathcal{F}) \rightarrow H^p(X,\mathcal{F}) $$
\emph{induced} by the couple $ (g,u) $ is defined as the composition of the maps $ H^p(g) $ and $ H^p(u) $.
\end{dfn}

\begin{rmk}
In case $ \mathcal{F} = \mathcal{O}_X $, there is a canonical isomorphism $ g^* \mathcal{O}_X \cong \mathcal{O}_X $, associated to the morphism $g$. So we get naturally an endomorphism of the cohomology groups
 $$ H^p(g) : H^p(X,\mathcal{O}_X) \rightarrow H^p(X,\mathcal{O}_X), $$
for all $p \geq 0$.
\end{rmk}

\subsection{}\label{5.4}
Let $ G $ be a finite group acting on $X$, and let $ \mathcal{F} $ be a coherent $ \mathcal{O}_X $-module. An isomorphism $ u : g^* \mathcal{F} \to \mathcal{F} $ is called a \emph{covering homomorphism}. We say that $ \mathcal{F} $ is a \emph{$G$-sheaf} if there exist, for every $ g \in G $, covering homomorphisms $u_g$ such that $ u_h \circ h^* u_g = u_{gh} $, where $ g, h \in G $.

\begin{rmk}
Let $G$ be a finite group acting on a projective scheme $X/\Spec(k)$, and let $ \mathcal{F} $ be a $G$-sheaf on $X$. The compatibility conditions ensure that $H^p(X,\mathcal{F})$ is a $ k[G] $-module for all $ p \geq 0 $.
\end{rmk}

Let $ u' : g^* \mathcal{F}' \to \mathcal{F}' $ be a second covering homomorphism. A map of covering homomorphisms is a map $ \phi : \mathcal{F} \to \mathcal{F}' $ of $ \mathcal{O}_X $-modules such that $ \phi \circ u = u' \circ g^* \phi $. A map of $G$-sheaves $ \mathcal{F} $ and $ \mathcal{F}' $ is a map of $ \mathcal{O}_X $-modules respecting the respective $G$-sheaf structures. The category of $G$-sheaves on $X$ is in fact an abelian category. For this, and further properties of $G$-sheaves, we refer to \cite{Kock}, Chapter 1.

Consider a short exact sequence 
$$ 0 \to (\mathcal{F}_1,u_1) \to (\mathcal{F}_2,u_2) \to (\mathcal{F}_3,u_3) \to 0 $$
of coverings. It is straight forward to check that this gives a commutative diagram
\begin{equation} 
\xymatrix{
\ldots \ar[r] & H^p(X, \mathcal{F}_2) \ar[r] \ar[d]^{H^p(g,u_2)} & H^p(X, \mathcal{F}_3) \ar[r]^{\delta} \ar[d]^{H^p(g,u_3)} & H^{p+1}(X, \mathcal{F}_1) \ar[r] \ar[d]^{H^{p+1}(g,u_1)} & \ldots \\
\ldots \ar[r] & H^p(X, \mathcal{F}_2) \ar[r] & H^p(X, \mathcal{F}_3) \ar[r]^{\delta} & H^{p+1}(X, \mathcal{F}_1) \ar[r] & \ldots .}
\end{equation}
Similarly, if $G$ acts on the projective scheme $X/k$, one checks that a short exact sequence of $G$-sheaves gives a long exact sequence of $k[G]$-modules in cohomology.




\subsection{}
If $X$ is projective over a field $k$, and $ \mathcal{F} $ is a coherent $ \mathcal{O}_X $-module, the cohomology groups $ H^p(X,\mathcal{F}) $ are finite dimensional $k$-vector spaces, and the trace $ \Tr(H^p(g,u)) $ of the endomorphism $ H^p(g,u) $ is defined. 
If $ \mathcal{F} $ is a $G$-sheaf, we let $ [H^p(X,\mathcal{F})] $ denote the element associated to $ H^p(X,\mathcal{F}) $ in the representation ring $R_G(k)$.

\subsection{}


In the case where $ p = \Char(k) > 0 $, we let $W(k)$ denote the ring of \emph{Witt vectors} for $k$ (\cite{Serre}, Chap.~II, par.~5). Recall that $W(k)$ is a complete discrete valuation ring, $p$ is a uniformizing parameter in $W(k)$ and the residue field is $k$. The fraction field $ FW(k) $ of $W(k)$, however, has characteristic $0$.

There exists a unique multiplicative map $ w : k \rightarrow W(k) $ that sections the reduction map $ W(k) \rightarrow k $. The map $w$ is often referred to as the \emph{Teichm\"uller lifting} from $k$ to $ W(k) $. 

Since we assume $ k = \overline{k} $, it follows that $k$ has a full set of $n$-th roots of unity, for any $n$ not divisible by $p$. As $W(k)$ is complete, these lift uniquely to $W(k)$, and reduction modulo $p$ induces an isomorphism of $ \boldsymbol{\mu}_n(W(k)) $ onto $ \boldsymbol{\mu}_n(k) $. 

\subsection{}



A few facts regarding \emph{Brauer characters} are needed, and are stated here in the case where $ G = \boldsymbol{\mu}_n $. We refer to \cite{SerreLin}, Chap.~18 for details.  

If $ E $ is a $ k[G] $-module, we let $ g_E $ denote the endomorphism of $ E $ induced by $ g \in G $. Since the order of $g$ divides $n$, and $n$ is relatively prime to $p$, it follows that $g_E$ is diagonalizable, and that all the eigenvalues $ \lambda_1, \ldots, \lambda_{e=\Dim E} $ are $n$-th roots of unity. The \emph{Brauer character} is then defined by assigning 
$$ \phi_E(g) = \sum_{i=1}^{e} w(\lambda_i). $$ 
It can be seen that the function $ \phi_E : G \rightarrow W(k) $ thus obtained is a class function on $G$. We shall call the element $ \phi_E(g) \in W(k) $ the \emph{Brauer trace} of $g_E$. The ordinary trace is obtained from the Brauer trace by reduction modulo $p$. 

An important property of the Brauer character is that it is additive on short exact sequences. That is, if 
$ E' \rightarrow E \rightarrow E''  $ is a short exact sequence of $ k[G] $-modules, then $ \phi_E = \phi_{E'} + \phi_{E''} $. A useful consequence is that if 
$$ 0 \rightarrow E_0 \rightarrow \ldots \rightarrow E_i \rightarrow \ldots \rightarrow E_l \rightarrow 0 $$
is an exact sequence of $k[G]$-modules, we get that $ \sum_{i=0}^l (-1)^i \phi_{E_i}(g) = 0 $.

\begin{ntn}
If $ V $ is a finite dimensional vector space over $k$, and $ \psi : V \rightarrow V $ is an automorphism, we will use the notation $ \Tr_{\beta}(\psi) $ for the Brauer trace of $\psi$.
\end{ntn}

\subsection{}
We now consider a smooth, connected and projective curve $C$ over $k$, with an invertible sheaf $ \mathcal{L} $ on $C$. 

Let $ g : C \rightarrow C $ be an automorphism, and let $ u : g^*\mathcal{L} \rightarrow \mathcal{L} $ be a covering map. We would like to compute the alternating sum $ \sum_{p=0}^1 (-1)^p ~\Tr_{\beta}(H^p(g,u)) $ of the Brauer traces. 

\subsection{}
Let us first consider the case when the automorphism $ g : C \rightarrow C $ is trivial, i.e., $ g = \Id_C $. Then $ H^p(g) $ is the identity, so we need only consider $ u : \mathcal{L} = g^*\mathcal{L} \rightarrow \mathcal{L} $. Hence $ H^p(g,u) = H^p(u) : H^p(C,\mathcal{L}) \rightarrow H^p(C,\mathcal{L}) $, where $ u \in \Aut_{\mathcal{O}_C}(\mathcal{L}) $. Since $ \Aut_{\mathcal{O}_C}(\mathcal{L}) = k^* $, we have that $u$ is multiplication with some element $ \lambda_u \in k^* $. 

\begin{prop}\label{prop. 5.5}
Let us keep the hypotheses above. Then the following equality holds in $W(k)$:
$$ \sum_{p=0}^1 (-1)^p~\Tr_{\beta}(H^{p}(u)) = w(\lambda_u) \cdot (\Deg_C(\mathcal{L}) + 1 - p_a(C)). $$
\end{prop}
\begin{pf}
Using \v{C}ech cohomology, it is straightforward to see that $H^p(u)$ is multiplication by $ \lambda_u $, for all $p$. Applying the Riemann-Roch formula then gives the result.
\end{pf}

\subsection{Computing the trace when $g$ is non-trivial}
Let now $ g \in \Aut_k(C) $ be a non-trivial automorphism of finite period $n$ (i.e. $ g^n = \Id_C $), where $n$ is not divisible by the characteristic of $k$. 


In this situation, the so called \emph{Lefschetz-Riemann-Roch} formula (\cite{Don}, Theorem 5.4, Corollary 5.5) gives a formula for the Brauer trace of $ H^p(g,u) $ in terms of local data at the fixed points of $g$. 

Let $ z \in C $ be a fixed point, and let $ i_z : \{z\} \hookrightarrow C $ be the inclusion. Pulling back $ u $ via $ i_z $ gives
$$ u(z) = i_z^* u : i_z^* g^* \mathcal{L} = i_z^* \mathcal{L} \rightarrow i_z^* \mathcal{L}, $$ 
a $k$-linear endomorphism of $ \mathcal{L}(z) $. We let $ \lambda_u(z) $ denote the (unique) eigenvalue of $ u(z) $

Since $ z \in C $ is a fixed point, there is an induced automorphism 
$$ dg(z) : \mathfrak{m}_z/\mathfrak{m}_z^2 \rightarrow \mathfrak{m}_z/\mathfrak{m}_z^2 $$ 
of the cotangent space of $C$ at $z$. We let $ \lambda_{dg}(z) $ denote the (unique) eigenvalue of $ dg(z) $.

The Lefschetz-Riemann-Roch then comes out as follows: 

\begin{prop}\label{prop. 5.6}
Let $C$, $ \mathcal{L} $, $g$ and $u$ be as above, and denote by $ C^g $ the (finite) set of fixed points of $g$. Then the following equality holds in $W(k)$: 
$$ \sum_{p=0}^1 (-1)^p ~\Tr_{\beta}(H^p(g,u)) = \sum_{z \in C^g} w(\lambda_u(z))/(1- w(\lambda_{dg}(z))). $$
\end{prop}
\begin{pf}
See Prop.~5.8 in~\cite{Thesis}.
\end{pf}

\begin{rmk}
The reader might want to compare Proposition \ref{prop. 5.6} with the \emph{Woods-Hole}-formula (\cite{SGA5}, Exp. III, Cor. 6.12), that gives a formula for the ordinary trace, instead of the Brauer trace, but with weaker assumptions on the automorphism $g$. 
\end{rmk}

\begin{rmk} 
Throughout the rest of the text we will, when no confusion can arise, continue to write $ \lambda $ instead of $w(\lambda) $ for the Teichm\"uller lift of a root of unity $\lambda$.
\end{rmk}

\section{Action on the minimal desingularization}\label{section 6}
Recall the set-up in Section \ref{extensions and actions}. We considered an SNC-model $ \mathcal{X}/S $, and a tamely ramified extension $ S'/S $ of degree $n$ that is prime to the least common multiple of the multiplicities of the irreducible components of $ \mathcal{X}_k $. The minimal desingularization of the pullback $ \mathcal{X}_{S'}/S' $ is an SNC-model $ \mathcal{Y}/S' $, and the Galois group $ G = \boldsymbol{\mu_n} $ of the extension $ S'/S $ acts on $ \mathcal{Y} $. 



Our goal is to compute the irreducible characters for this representation on $ H^1(\mathcal{Y}_k, \mathcal{O}_{\mathcal{Y}_k}) $. To do this, we would ideally compute the Brauer trace of the automorphism of $ H^1(\mathcal{Y}_k, \mathcal{O}_{\mathcal{Y}_k}) $ induced by $g$, for every group element $ g \in G $. This information would then be used to compute the Brauer character. However, we can not do this directly. Instead we will compute the Brauer trace of the automorphism induced by $g$ on the formal difference $ H^0(\mathcal{Y}_k, \mathcal{O}_{\mathcal{Y}_k}) - H^1(\mathcal{Y}_k, \mathcal{O}_{\mathcal{Y}_k}) $, for any $ g \in G $. In our applications, we know the character for $ H^0(\mathcal{Y}_k, \mathcal{O}_{\mathcal{Y}_k}) $, so this would suffice in order to determine the character for $ H^1(\mathcal{Y}_k, \mathcal{O}_{\mathcal{Y}_k}) $.

The fact that $ \mathcal{Y}_k $ is not in general smooth, prevents us from using Propositions \ref{prop. 5.5} and \ref{prop. 5.6} directly. On the other hand, the irreducible components of $ \mathcal{Y}_k $ are smooth and proper curves. So we shall in fact show that it is possible to reduce to computing Brauer traces on each individual component of $ \mathcal{Y}_k $, where Propositions \ref{prop. 5.5} and \ref{prop. 5.6} do apply. The key step in obtaining this is to introduce a certain filtration of the special fiber $\mathcal{Y}_k$.

\subsection{} 
Let $ \{ C_{\alpha} \}_{\alpha \in \mathcal{A}} $ denote the set of irreducible components of $ \mathcal{Y}_k $, and let $ m_{\alpha} $ denote the multiplicity of $ C_{\alpha} $ in $ \mathcal{Y}_k $. Then $ \mathcal{Y}_k $ can be written in Weil divisor form as
$$ \mathcal{Y}_k = \sum_{\alpha} m_{\alpha} C_{\alpha}. $$

\begin{dfn}\label{complete}
A \emph{complete} filtration of $\mathcal{Y}_k$ is a sequence 
$$ 0 < Z_m < \ldots < Z_j < \ldots < Z_1 = \mathcal{Y}_k $$
of effective divisors $Z_j$ supported on $\mathcal{Y}_k$, such that for each $1 \leq j \leq m - 1 $ there exists an $ \alpha_j \in \mathcal{A} $ with $ Z_j - Z_{j+1} = C_{\alpha_j} $. So $ m = \sum_{\alpha} m_{\alpha} $.
\end{dfn}

Loosely speaking, such a filtration of $\mathcal{Y}_k$ is obtained by removing the irreducible components of the special fiber one at the time (counted with multiplicity).

\subsection{}
At each step $Z' < Z$ of a complete filtration, we can construct an exact sequence of sheaves. 

\begin{lemma}\label{lemma 6.1}
Let $ 0 \leq Z' < Z \leq \mathcal{Y}_k $ be divisors such that $ Z - Z' = C $, for some irreducible component $C$ of $ \mathcal{Y}_k $. Denote by $ \mathcal{I}_Z $ and $  \mathcal{I}_{Z'} $ the corresponding ideal sheaves in $ \mathcal{O}_{\mathcal{Y}} $. Let $ i_{Z} $, $ i_{Z'} $ and $ i_C $ be the canonical inclusions of $ Z $, $Z'$ and $ C $ in $\mathcal{Y}$. Furthermore, let $ \mathcal{L} = i_C^*(\mathcal{I}_{Z'}) $. We then have an exact sequence 
$$ 0  \rightarrow (i_C)_* \mathcal{L} \rightarrow (i_{Z})_* \mathcal{O}_{Z} \rightarrow (i_{Z'})_* \mathcal{O}_{Z'} \rightarrow 0 $$
of $ \mathcal{O}_{\mathcal{Y}} $-modules. 
\end{lemma}
\begin{pf}
The inclusions $ \mathcal{I}_Z \subset \mathcal{I}_{Z'} \subset \mathcal{O}_{\mathcal{Y}} $ give rise to an exact sequence
$$ 0 \rightarrow \mathcal{K} \rightarrow \mathcal{O}_{\mathcal{Y}}/\mathcal{I}_Z \rightarrow \mathcal{O}_{\mathcal{Y}}/\mathcal{I}_{Z'} \rightarrow 0, $$
where $ \mathcal{K} = \mathcal{I}_{Z'} / \mathcal{I}_Z $ denotes the kernel. We need to determine $ \mathcal{K} $. 

Consider the surjection $ \mathcal{I}_{Z'} \rightarrow \mathcal{I}_{Z'} / \mathcal{I}_Z $. Pulling back with $ i_C^* $, we get a surjection $ i_C^*(\mathcal{I}_{Z'}) \rightarrow i_C^* (\mathcal{I}_{Z'} / \mathcal{I}_Z) $, and we claim that this map is an isomorphism. Indeed, let $ U = \Spec(A) \subset \mathcal{Y} $ be an open affine set. Then $A$ is a regular domain, and the ideal sheaves $ \mathcal{I}_C $, $ \mathcal{I}_Z $ and $ \mathcal{I}_{Z'} $ restricted to $U$ correspond to invertible modules $ I_C $, $ I_{Z} $ and $ I_{Z'} $ in $A$. Since $ Z = Z' + C $, we have that $ I_{Z} = I_C I_{Z'} $. So $ I_{Z'}/I_{Z} = I_{Z'}/I_C I_{Z'} = I_{Z'} \otimes_A A/I_C $. From this observation, it follows easily that the map above is an isomorphism on all stalks, and therefore an isomorphism.
\end{pf} 

\subsection{}\label{6.4} The filtrations we have introduced are $G$-equvariant.


\begin{lemma}\label{G-sheaf} 
Let $ 0 \leq Z \leq \mathcal{Y}_k $ be an effective divisor, with ideal sheaf $ \mathcal{I}_Z $. Then we have that $ \mathcal{I}_Z $ is a $G$-subsheaf of $ \mathcal{O}_{\mathcal{Y}} $.
\end{lemma}
\begin{pf}
Let $ g \in G $ be any group element. Applying the exact functor $ g^{-1} $ to the inclusion $ \mathcal{I}_Z \subset \mathcal{O}_{\mathcal{Y}} $ gives an inclusion $ g^{-1} \mathcal{I}_Z \subset g^{-1} \mathcal{O}_{\mathcal{Y}} $. Composing this inclusion with the canonical map $ g^{\sharp} : g^{-1} \mathcal{O}_{\mathcal{Y}} \rightarrow \mathcal{O}_{\mathcal{Y}} $, we obtain a map $g^{-1} \mathcal{I}_Z \rightarrow \mathcal{O}_{\mathcal{Y}}$. 
Now, let $ \mathcal{J} $ be the sheaf of ideals generated by the image of $ g^{-1} \mathcal{I}_Z $ in $ \mathcal{O}_{\mathcal{Y}} $. We have that $ \mathcal{J} $ is the ideal sheaf of $ g^{-1}(Z) $. But in our case $ g^{-1}(Z) = Z $, and therefore $ \mathcal{J} = \mathcal{I}_Z $.

The inclusion above induces an injective map $ g^*\mathcal{I}_Z \rightarrow g^* \mathcal{O}_{\mathcal{Y}} \cong \mathcal{O}_{\mathcal{Y}} $ of $ \mathcal{O}_{\mathcal{Y}} $-modules, whose image is $ \mathcal{I}_Z = g^{-1} \mathcal{I}_Z \cdot \mathcal{O}_{\mathcal{Y}} $ (\cite{Hart}, II.7.12.2). Hence we obtain an isomorphism $ u_Z : g^*\mathcal{I}_Z \rightarrow \mathcal{I}_Z $.

It is easy to check that the isomorphisms $ g^*\mathcal{I}_Z \rightarrow \mathcal{I}_Z $ for various elements $ g \in G $ satisfy the compatibility conditions, and commute with the $G$-sheaf structure on $ \mathcal{O}_{\mathcal{Y}} $.
\end{pf}

\subsection{}\label{6.7}

\begin{prop}\label{prop. 6.5}
Let us keep the hypotheses and notation from Lemma \ref{lemma 6.1}. The sequence 
$$ 0  \rightarrow (i_C)_* \mathcal{L} \rightarrow (i_{Z})_* \mathcal{O}_{Z} \rightarrow (i_{Z'})_* \mathcal{O}_{Z'} \rightarrow 0 $$
is an exact sequence of $ G $-sheaves.
\end{prop}
\begin{pf}
Let $ \mathcal{I}_{Z} \subset \mathcal{I}_{Z'} \subset \mathcal{O}_{\mathcal{Y}} $ be the inclusions of the ideal sheaves. From Lemma \ref{G-sheaf}, it follows that these maps are maps of $G$-sheaves. The result now follows from the fact that the category of $G$-modules on $\mathcal{Y}$ is an abelian category (\cite{Kock}, Lemma 1.3). 
\end{pf}

\vspace{0.3cm}

In particular, Proposition \ref{prop. 6.5} implies that for any $g \in G$, there are covering maps $u$ (resp.~$v$, $v'$) of $ (i_C)_* \mathcal{L} $ (resp.~$(i_Z)_* \mathcal{O}_{Z}$, $(i_{Z'})_* \mathcal{O}_{Z'}$) giving an exact sequence
\begin{equation}\label{G-shortseq}
0 \to ((i_C)_* \mathcal{L}, u) \to ((i_Z)_* \mathcal{O}_{Z},v) \to ((i_{Z'})_* \mathcal{O}_{Z'},v') \to 0 
\end{equation}
of covering maps.



The maps $ u $, $ v $ and $ v' $ induce, for every $ p \geq 0 $, automorphisms $ H^p(g,u) $, $ H^p(g,v) $ and $ H^p(g,v') $ that commute with the differentials in the long exact sequence in cohomology asociated to Sequence (\ref{G-shortseq}). That is, we obtain, for every $ g \in G $, an \emph{automorphism} of this long exact sequence. 

Note that all the sheaves appearing in the exact sequence above are supported on the special fiber of $ \mathcal{Y} $. We will now explain how we can ``restrict'' the endomorphisms $ H^p(g,u) $, $ H^p(g,v) $ and $ H^p(g,v') $ to the support of the various sheaves.


\subsection{}\label{6.9}

Let $X$ and $Y$ be schemes, and let $ i : X \rightarrow Y $ be a closed immersion. Assume also that an automorphism $ g : Y \rightarrow Y $ is given, that restricts to an automorphism $ f = g|_X : X \rightarrow X $. 

Note that if $ \mathcal{F} $ is a quasi-coherent sheaf on $X$, then the push-forward $ i_* \mathcal{F} $ is a quasi-coherent sheaf on $ Y $, since $i$ is a closed immersion. The following lemma is straightforward, yet tedious to prove, so we omit the proof. 

\begin{lemma}\label{lemma 6.6}
Keep the hypotheses above. Let $ \mathcal{F} $ be a quasi-coherent sheaf on $X$, and let $ u : g^* i_* \mathcal{F} \rightarrow i_* \mathcal{F} $ be a homomorphism of $ \mathcal{O}_Y $-modules. Then there is induced, for every $ p \geq 0 $, a commutative diagram
$$ \xymatrix{
H^p(Y, i_* \mathcal{F}) \ar[rr]^{H^p(g,u)} \ar[d]_{H^p(i)} & & H^p(Y, i_* \mathcal{F}) \ar[d]^{H^p(i)} \\
H^p(X, i^* i_*\mathcal{F}) \ar[rr]^{H^p(f,i^*u)} & & H^p(X, i^* i_*\mathcal{F}),} 
$$
where the vertical arrows are isomorphisms.
\end{lemma}
\begin{pf}
See Prop. 6.7 in \cite{Thesis}.
\end{pf}

\subsection{}\label{6.11}
The closed immersions $ i_C $, $ i_Z $ and $ i_{Z'} $ induce isomorphisms
$$ H^p(\mathcal{Y}, (i_C)_* \mathcal{L}) \cong H^p(C,\mathcal{L}), $$ 
$$ H^p(\mathcal{Y}, (i_Z)_* \mathcal{O}_Z) \cong H^p(Z, \mathcal{O}_Z) $$
and
$$ H^p(\mathcal{Y}, (i_{Z'})_* \mathcal{O}_{Z'}) \cong H^p(Z', \mathcal{O}_{Z'}), $$
for all $ p \geq 0 $. Here we have identified $ \mathcal{L} $ with $ (i_C)^* (i_C)_* \mathcal{L} $ (and likewise for $ \mathcal{O}_Z $ and $ \mathcal{O}_{Z'} $).

Since $ C $, $ Z $ and $ Z' $ are projective curves over $k$, and since $ \mathcal{L} $, $ \mathcal{O}_Z $ and $ \mathcal{O}_{Z'} $ are coherent sheaves, the cohomology groups above are finite dimensional $k$-vector spaces, and nonzero only for $ p = 0 $ and $ p = 1 $.
So the long exact sequence in cohomology associated to Sequence (\ref{G-shortseq}) is simply
\begin{equation}\label{sequence 6.7} 
0 \rightarrow H^0(C,\mathcal{L}) \rightarrow H^0(Z, \mathcal{O}_Z) \rightarrow \ldots \to H^1(Z, \mathcal{O}_Z) \rightarrow H^1(Z', \mathcal{O}_{Z'}) \rightarrow 0.
\end{equation}

\begin{prop}\label{longexactgroup}
Sequence (\ref{sequence 6.7}) is an exact sequence of $k[G]$-modules. Furthermore, we get an equality
\begin{equation}
\sum_{p=0}^1 (-1)^p [H^p(Z, \mathcal{O}_{Z})] = \sum_{p=0}^1 (-1)^p [H^p(Z', \mathcal{O}_{Z'})] + \sum_{p=0}^1 (-1)^p [H^p(C,\mathcal{L})]  
\end{equation}
of (virtual) $k[G]$-modules.
\end{prop}
\begin{pf}
Denote by $g_C$ the restriction of $g$ to $C$. By Lemma \ref{lemma 6.6}, restriction to $C$ gives a commutative diagram
$$ \xymatrix{ H^p(\mathcal{Y},(i_C)_* \mathcal{L}) \ar[d]_{H^p(g,u)} \ar[rr]^{\cong} & & H^p(C,\mathcal{L}) \ar[d]^{H^p(g_C, i_C^*u)} \\
H^p(\mathcal{Y},(i_C)_* \mathcal{L}) \ar[rr]^{\cong} & & H^p(C,\mathcal{L}). } $$
Also, we get similar diagrams for $ H^p(g_Z, i_Z^* v) $ and $ H^p(g_{Z'}, i_{Z'}^* v') $. 

Having made these identifications, we see that the automorphisms $ H^p(g_C, i_C^* u)$, $ H^p(g_Z, i_Z^* v) $ and $ H^p(g_{Z'}, i_{Z'}^* v') $, for $ p = 0, 1 $, fit together to give an automorphism of Sequence \ref{sequence 6.7} above. One checks that Sequence (\ref{sequence 6.7}) is an exact sequence of $ k[G] $-modules, and from this fact, the second statement immediately follows.
\end{pf}


\subsection{}
Let us write $ \mathcal{Y}_k = \sum_{\alpha} m_{\alpha} C_{\alpha} $, where $ \alpha \in \mathcal{A} $, and put $ m = \sum_{\alpha} m_{\alpha} $. Fix a complete filtration
$$ 0 < Z_m < \ldots < Z_j < \ldots < Z_2 < Z_1 = \mathcal{Y}_k, $$ 
where $ Z_j - Z_{j+1} = C_j $ for some $ C_j \in \{ C_{\alpha} \}_{\alpha \in \mathcal{A} } $, for each $ j \in \{ 1, \ldots, m-1 \} $. At each step of this filtration, Lemma \ref{lemma 6.1} asserts that there is a short exact sequence
$$ 0 \rightarrow (i_{C_j})_* \mathcal{L}_j \rightarrow (i_{Z_j})_* \mathcal{O}_{Z_j} \rightarrow (i_{Z_{j+1}})_* \mathcal{O}_{Z_{j+1}} \rightarrow 0, $$
where $ i_{\star} : \star \hookrightarrow \mathcal{Y} $ is the canonical inclusion. Note in particular that $ Z_m = C_m $, for some $ C_m \in \{ C_{\alpha} \}_{\alpha \in \mathcal{A} } $, so it makes sense to  write $ \mathcal{O}_{Z_m} = \mathcal{L}_m $.


Proposition \ref{longexactgroup} has the following nice consequence: 

\begin{prop}\label{prop-Grep}
We have an equality
\begin{equation}
\sum_{p=0}^1 (-1)^p [H^p(\mathcal{Y}_k, \mathcal{O}_{\mathcal{Y}_k})] = \sum_{j=1}^m \sum_{p=0}^1 (-1)^p [H^p(C_j, \mathcal{L}_j)] 
\end{equation}
of $ k[G] $-modules.
\end{prop}
\begin{pf}
Follows easily by induction from Proposition \ref{longexactgroup}.
\end{pf}

\vspace{0.3 cm}

Let $ \phi^p_g $ (resp.~$\phi^p_{g,j}$) be the automorphism induced by the element $ g \in G $ on $ H^p(\mathcal{Y}_k, \mathcal{O}_{\mathcal{Y}_k}) $ (resp.~$H^p(C_j, \mathcal{L}_j)$). The following corollary is immediate: 

\begin{cor}\label{cor-Grep}
For any $g \in G$, the following equality holds in $W(k)$:
$$ \sum_{p=0}^1 (-1)^p~\Tr_{\beta}(\phi^p_g) = \sum_{j=1}^m \sum_{p=0}^1 (-1)^p~\Tr_{\beta}(\phi^p_{g,j}). $$
\end{cor}




\subsection{}
The importance of Corollary \ref{cor-Grep} is that it reduces the problem of computing the alternating sum of the Brauer traces of the endomorphisms 
$$ \phi^p_g : H^p(\mathcal{Y}_k, \mathcal{O}_{\mathcal{Y}_k}) \to H^p(\mathcal{Y}_k, \mathcal{O}_{\mathcal{Y}_k}) $$  
to instead computing the same data for the endomorphisms 
$$ \phi^p_{g,j} : H^p(C_j, \mathcal{L}_j) \to H^p(C_j, \mathcal{L}_j) $$
for certain invertible sheaves $ \mathcal{L}_j $, supported on the smooth irreducible components $ C_j $ of $ \mathcal{Y}_k $. 
 The main benefit is that for the latter computations, we can apply the Lefschetz-Riemann-Roch formulas in Proposition \ref{prop. 5.5} and \ref{prop. 5.6}. In what follows, we will explain how this can be done.

We keep the notation from Lemma \ref{lemma 6.1}, hence we have $ 0 \leq Z' < Z \leq \mathcal{Y}_k $, where $ Z - Z' = C $. We let $ g \in G $ be an element corresponding to a primitive root of $ \boldsymbol{\mu}_n $. From Proposition \ref{prop. 3.3}, we see that the fixed points of the automorphism $ g : C \to C $ are precisely the two points where $ C $ meets the other components of $ \mathcal{Y}_k $. Let $ y \in C $ be one of the fixed points, and let $ dg_y $ denote the cotangent map at $y$. The eigenvalue of $ dg_y $ can easily be computed using Proposition \ref{prop. 3.3}. 

We will also need to compute the eigenvalue of the induced automorphism $ u_y : \mathcal{L}(y) \to \mathcal{L}(y) $. To do this, let $ C' $ be the other component of $ \mathcal{Y}_k $ that passes through $ y $. Then we can write $ \mathcal{I}_{Z'} = \mathcal{I}_{C}^{\otimes a} \otimes \mathcal{I}_{C'}^{\otimes a'} \otimes \mathcal{I}_{0} $, where $ \mathcal{I}_{0} $ is the ideal sheaf of an effective Cartier divisor not containing $C$ or $C'$. 

Since $C$ and $C'$ intersect transversally at $y$, the fibers $ \mathcal{I}_{C}(y) $ and $ \mathcal{I}_{C'}(y) $ generate the cotangent space to $ \mathcal{Y} $ at $y$. The eigenvalues $ \lambda $ and $ \lambda' $ of these generators can easily be computed using Proposition \ref{prop. 3.3}, since they correspond to the two coordinates locally at $y$. The following lemma is an easy computation (c.f.~\cite{Thesis}, Section 6.10):

\begin{lemma}\label{eigenvaluelemma}
Keep the notation from the discussion above. The unique eigenvalue of the automorphism $ u_y : \mathcal{L}(y) \to \mathcal{L}(y) $ is $ \lambda^a \lambda'^{a'} $.
\end{lemma}

\section{Special filtrations for trace computations}\label{special filtrations}

Let $ \mathcal{X}/S $ be an SNC-model, and let $ S' \rightarrow S $ be a tame extension of degree $n$, where $n$ is prime to the least common multiple of the multiplicities of the irreducible components of $ \mathcal{X}_k $. Let $ \mathcal{X}' $ be the normalization of $ \mathcal{X}_{S'} $, and let $ \mathcal{Y} $ be the minimal desingularization of $ \mathcal{X}' $. 

This section is devoted to computing, for any $ g \in G = \boldsymbol{\mu}_n $, the Brauer trace of the automorphism induced by $g$ on the formal difference $ H^0(\mathcal{Y}_k, \mathcal{O}_{\mathcal{Y}_k}) - H^1(\mathcal{Y}_k, \mathcal{O}_{\mathcal{Y}_k}) $. Hence a lot of our previous work will come together in this section. 

Our assumption on the degree of $ S'/S $ makes it possible to describe $ \mathcal{Y}_k $ in terms of $ \mathcal{X}_k $. In particular, since every component of $ \mathcal{Y}_k $ either is an exceptional curve, or dominates a component of $ \mathcal{X}_k $, it is natural to stratify the combinatorial structure of $ \mathcal{Y}_k $ according to the combinatorial structure of $ \mathcal{X}_k $.

This stratification proves to be very convenient for our trace computations. The section concludes with Theorem \ref{thm. 9.13}, which gives a formula for the trace mentioned above as a sum of contributions associated in a natural way to the combinatorial structure of $ \mathcal{X}_k $. 

\subsection{}\label{9.1}
We will associate a graph $ \Gamma(\mathcal{X}_k) $ to $ \mathcal{X}_k $ in the following way: The set of vertices, $ \mathcal{V} $, consists of the irreducible components of $ \mathcal{X}_k $. The set of edges, $ \mathcal{E} $, consists of the intersection points of $ \mathcal{X}_k $, and two distinct vertices $\upsilon$ and $\upsilon'$ are connected by $ \Card (\{ D_{\upsilon} \cap D_{\upsilon'} \}) $ edges, where $D_{\upsilon}$ denotes the irreducible component corresponding to $\upsilon$.

We define two natural functions on the set of vertices $ \mathcal{V} $. First, let the \emph{genus}
$$ \mathfrak{g} : \mathcal{V} \rightarrow \mathbb{N}_0, $$ 
be defined by $ \mathfrak{g}(\upsilon) = p_a(D_{\upsilon}) $. We also let the \emph{multiplicity} 
$$ \mathfrak{m} : \mathcal{V} \rightarrow \mathbb{N}, $$ 
be defined by $ \mathfrak{m}(\upsilon) = \Mult_{\mathcal{X}_k}(D_{\upsilon}) $.

The graph $ \Gamma(\mathcal{X}_k) $, together with the functions $\mathfrak{g}$ and $\mathfrak{m}$, encode all the combinatorial and numerical properties of $ \mathcal{X}_k $.

\subsection{}\label{9.2}

Let $ \mathcal{S} $ denote the set of irreducible components of $ \mathcal{Y}_k $. If $ C \in \mathcal{S} $, then we have either: 
\begin{enumerate}
\item $ C $ dominates a component $ D_{\upsilon} $ of $ \mathcal{X}_k $, or  

\item $ C $ is a component of the exceptional locus of the minimal desingularization $ \rho : \mathcal{Y} \rightarrow \mathcal{X}' $.
\end{enumerate}

In the first case, we have that $ p_a(C) = \mathfrak{g}(\upsilon) $, and $ \Mult_{\mathcal{Y}_k}(C) = \mathfrak{m}(\upsilon) $. Furthermore, $G$ acts trivially on $C$. Since $C$ is the unique component of $ \mathcal{Y}_k $ corresponding to $ \upsilon $, we write $C = C_{\upsilon}$.

In the second case, we have that $ C $ is part of a chain of exceptional curves, corresponding uniquely to an edge $ \varepsilon \in \mathcal{E} $. Hence $ p_a(C) = 0 $. By choosing an ordering (or direction) of this chain, we can index the components in the chain by $ l $, for $ 1 \leq l \leq L(\varepsilon) $, where $ L(\varepsilon) $ is the length of the chain. So we can write $ C = C_{\varepsilon,l} $, for some $l \in \{1, \ldots, L(\varepsilon) \} $. By Proposition \ref{prop. 3.3}, $G$ acts nontrivially on $ C $, with fixed points exactly at the two points where $C$ meets the rest of the special fiber.

The special fiber $ \mathcal{Y}_k $ can now be written, as an effective divisor on $ \mathcal{Y} $, in the form 
$$ \mathcal{Y}_k = \sum_{\varepsilon \in \mathcal{E}} \sum_{l = 1}^{L(\varepsilon)} \mu_{\varepsilon,l} C_{\varepsilon,l} + \sum_{\upsilon \in \mathcal{V}} m_{\upsilon} C_{\upsilon}, $$
where $ \mu_{\varepsilon,l} $ denotes the multiplicity of the component $ C_{\varepsilon,l} $, and $m_{\upsilon}$ is the multiplicity of $C_{\upsilon}$.

\subsection{}\label{9.3}
We will now consider \emph{special} filtrations of $ \mathcal{Y}_k $, inspired by the partition of the set of irreducible components of $ \mathcal{Y}_k $ introduced above. 

Let us choose an ordering of the elements in $ \mathcal{V} $. We can then define the following sequence:
$$ 0 < \ldots < Z_{\mathcal{E}} =: Z_{\upsilon_{|\mathcal{V}| + 1}} < Z_{\upsilon_{|\mathcal{V}|}} < \ldots < Z_{\upsilon_i} < \ldots < Z_{\upsilon_1} = \mathcal{Y}_k, $$
where $ Z_{\mathcal{E}} := \mathcal{Y}_k - \sum_{\upsilon \in \mathcal{V}} m_{\upsilon} C_{\upsilon} $. The $ Z_{\upsilon_i} $ are defined inductively, for every $ i \in \{ 1, \ldots, |\mathcal{V}| \} $, by the refinements
$$ Z_{\upsilon_{i+1}} = Z_{\upsilon_i}^{m_{\upsilon_i} + 1} < \ldots < Z_{\upsilon_i}^j < \ldots < Z_{\upsilon_i}^1 = Z_{\upsilon_i}, $$ 
where $ Z_{\upsilon_i}^{j + 1} = Z_{\upsilon_i} - j C_{\upsilon_i} $ for every $ j \in \{ 0, \ldots, m_{\upsilon_i} \} $.

Next, we choose an ordering of the elements in $ \mathcal{E} $. We can then define the following sequence:
$$ 0 =: Z_{\varepsilon_{|\mathcal{E}|+1}} < Z_{\varepsilon_{|\mathcal{E}|}} < \ldots < Z_{\varepsilon_i} < \ldots < Z_{\varepsilon_1} := Z_{\mathcal{E}}. $$  
 The $ Z_{\varepsilon_i} $ are defined inductively, for any $ i \in \{ 1, \ldots, |\mathcal{E}| \} $, by the refinements
$$ Z_{\varepsilon_{i+1}} := Z_{\varepsilon_i, L(\varepsilon_i) + 1} < \ldots < Z_{\varepsilon_i,l} < \ldots < Z_{\varepsilon_i,1} := Z_{\varepsilon_i}, $$
which in turn are defined inductively, for every $ l \in \{ 1, \ldots, L(\varepsilon_i) \} $, by the further refinements
$$ Z_{\varepsilon_i, l+1} := Z_{\varepsilon_i, l}^{\mu_l + 1} < \ldots < Z_{\varepsilon_i, l}^j < \ldots < Z_{\varepsilon_i, l}^1 := Z_{\varepsilon_i, l}, $$
where $ Z_{\varepsilon_i, l}^{j+1} := Z_{\varepsilon_i, l} - j C_{\varepsilon_i, l} $, for every $ j \in \{ 0, \ldots, \mu_l \} $.

\begin{ex}\label{specfilex}
Let $ \mathcal{X}/S $ be an SNC-model with special fiber $ \mathcal{X}_k = 3 D_4 + D_1 + D_2 + D_3 $, where $D_i$ meets $D_4$ in a unique point for $ i \in \{ 1,2,3 \} $, and with no further intersection points. Let $R'/R$ be a tame extension of degree $7$, and let $ \mathcal{Y}/S' $ be the minimal desingularization of the normalization $ \mathcal{X}' $ of $ \mathcal{X}_{S'} $.

The singularities of $ \mathcal{X}' $ are formally isomorphic to $ \sigma = (1,3,7) $. From Example \ref{Gdesingex} we know that the exceptional locus of the resolution of $ \sigma $ consists of two components of multiplicity $1$. Let us write $C_i^1$ and $C_i^2$ for the components corresponding to the edge $ \varepsilon_i = (D_i,D_4) $. 

We can now write $ \mathcal{Y}_k = \sum_{i=1}^4 m_i C_i + \sum_{j=1}^3 (C_j^1 + C_j^2) $. The first part of a special filtration is then
$$ Z_{\mathcal{E}} := \mathcal{Y}_k - \sum_{i=1}^4 m_i C_i < \ldots < \mathcal{Y}_k - (C_1 + C_2) < \mathcal{Y}_k - C_1 < \mathcal{Y}_k, $$
and the second part looks like
$$ 0 < C_3^2 < C_3^1 + C_3^2 < \ldots < Z_{\mathcal{E}} - (C_1^1 + C_1^2) < Z_{\mathcal{E}} - C_1^1 < Z_{\mathcal{E}}. $$
\end{ex}

\subsection{} 
In the rest of this paper, we shall always choose complete filtrations of $ \mathcal{Y}_k $ that are of the form
\begin{equation}
0 < \ldots < Z_{\varepsilon_i} < \ldots < Z_{\varepsilon_1} = Z_{\mathcal{E}} < \ldots < Z_{\upsilon_i} < \ldots < Z_{\upsilon_1} = \mathcal{Y}_k, 
\end{equation}
where $ Z_{\upsilon_{i+1}} < Z_{\upsilon_i} $ and $ Z_{\varepsilon_{i+1}} < Z_{\varepsilon_i} $ are subfiltrations as described above. We shall soon see that the chosen orderings of the sets $ \mathcal{E} $ and $ \mathcal{V} $ are irrelevant. 

The nice feature of working with filtrations like this becomes evident when one wants to do trace computations \`a la Section \ref{section 6}. Then we may actually reduce to considering subfiltrations $ Z_{\upsilon_{i+1}} < Z_{\upsilon_i} $, which we interpret as \emph{contributions from the vertices} of $\Gamma$, and subfiltrations $ Z_{\varepsilon_{i+1}} < Z_{\varepsilon_i} $, which we interpret as \emph{contributions from the edges}.

\subsection{}
Let us fix a vertex $\upsilon \in \mathcal{V}$. We shall now define and calculate the contribution to the trace from $\upsilon$. To do this, we choose a filtration of $ \mathcal{Y}_k $ as in Section \ref{9.3} above. Then there will be a subfiltration of the form:
 $$ Z_{\mathcal{E}} \leq Z_{\upsilon}^{m_{\upsilon}+1} < \ldots < Z_{\upsilon}^k < \ldots < Z_{\upsilon}^1 = Z_{\upsilon} \leq \mathcal{Y}_k, $$
where $ Z_{\upsilon}^{k} - Z_{\upsilon}^{k + 1} = C_{\upsilon} $, for all $ 1 \leq k \leq m_{\upsilon} $. The invertible sheaf associated to the $k$-th step in this filtration is $ \mathcal{L}_{\upsilon}^k := j_{\upsilon}^*(\mathcal{I}_{Z_{\upsilon}^{k+1}}) $, where $ j_{\upsilon} : C_{\upsilon} \hookrightarrow \mathcal{Y} $ is the canonical inclusion.

We will use the following easy lemma, whose proof is omitted.
\begin{lemma}\label{lemma 9.1}
Assume that $ S'/S $ is a nontrivial extension. If $C_1$ and $C_2$ are two distinct components of $ \mathcal{Y}_k $, corresponding to elements in $ \mathcal{V} $, then they have empty intersection.
\end{lemma}

In what follows, we will suppress the index $ \upsilon $, to simplify notation. Let $ D_1, \ldots, D_f $ be the irreducible components of $ Z $ that intersect $ C $ non-trivially, and that are not equal to $ C $. Let $ a_i $ denote the multiplicity of $D_i$. It follows from Lemma \ref{lemma 9.1} that the $D_i$ are exceptional components. Moreover, it follows from the way we constructed the filtration that the $D_i$ are precisely the components of $ \mathcal{Y}_k $ different from $ C $ that have non-empty intersection with $ C $. We can then write 
$$ Z^{k + 1} = (m - k) C + a_1 D_1 + \ldots + a_f D_f + Z_0, $$
where all components of $ Z_0 $ have empty intersection with $ C $. So we get that
\begin{equation}\label{equation 9.2}
\mathcal{L}^k  =  j^*\mathcal{I}_{Z^{k+1}} = (\mathcal{I}_{C}|_{C})^{ \otimes m - k } \otimes (\mathcal{I}_{D_1}|_{C})^{\otimes a_1} \otimes \ldots \otimes (\mathcal{I}_{D_f}|_{C})^{\otimes a_f}. 
\end{equation} 

Let $ g $ be an element of $ G = \boldsymbol{\mu}_n $, corresponding to a root of unity $ \xi $. Note that the restriction of the automorphism $ g $ to $C$ is $ \Id_{C}$. Let 
$$ \phi^p_{g,k} : H^p(C,\mathcal{L}^k) \to H^p(C,\mathcal{L}^k) $$
be the automorphism induced by $g$.

\begin{dfn}\label{contribution from vertex}
We define \emph{the contribution to the trace} from the vertex $ \upsilon \in \mathcal{V} $ as the sum
$$ \Tr_{\upsilon}(\xi) = \sum_{k=1}^{m} \sum_{p=0}^1 (-1)^p~\Tr_{\beta}( \phi^p_{g,k} ). $$
\end{dfn}

\vspace{0.3 cm}

The proposition below gives an effective formula for the contribution from a vertex $ \upsilon $:

\begin{prop}\label{prop. 9.7} The contribution to the trace from the vertex $ \upsilon $ is given by the formula
$$ \Tr_{\upsilon}(\xi) = \sum_{k=0}^{m-1} (\xi^{\alpha_{m}})^{k} ((m - k) C^2 + 1 - p_a(C)), $$
where $ \alpha_{m} $ is an inverse to $ m $ modulo $n$. 
\end{prop}
\begin{pf}
By Proposition \ref{prop. 5.5}, we have that
$$ \sum_{p=0}^1 (-1)^p~\Tr_{\beta}( \phi^p_{g,k} ) = \lambda_k(\Deg_{C}(\mathcal{L}^k) + 1 - p_a(C)), $$
where $ \lambda_k $ is the eigenvalue of the automorphism $ \mathcal{L}(y) \to \mathcal{L}(y) $, for \emph{any} point $ y \in C $. The proof will consist of specifying precisely the terms appearing in this formula. 

Let us first compute $ \Deg_{C}(\mathcal{L}^k) $. Since $ \mathcal{I}_{C} = \mathcal{O}_{\mathcal{Y}}(- C) $, it follows that
$$ \Deg_{C}( \mathcal{I}_{C}|_{C}) = \Deg_{C}(\mathcal{O}_{\mathcal{Y}}(- C)|_{C}) = - \Deg_{C}(\mathcal{O}_{\mathcal{Y}}(C)|_{C}) = - C^2.$$
Furthermore, for any $ i \in \{ 1, \ldots, f \} $, we have that $ \mathcal{I}_{D_i} = \mathcal{O}_{\mathcal{Y}}(- D_i) $, and hence
$$ \Deg_{C}(\mathcal{I}_{D_i}|_{C}) = \Deg_{C}(\mathcal{O}_{\mathcal{Y}}(- D_i)|_{C}) = - 1. $$
It then follows from Equation \ref{equation 9.2} that
$$ \Deg_{C}(\mathcal{L}^k) = - (m - k) C^2 - (a_1 + \ldots + a_f). $$ 
On the other hand, we have that $ - C^2 = (a_1 + \ldots + a_f)/m $, and therefore we get that
$ \Deg_{C}(\mathcal{L}^k) = k C^2 $.

We now claim that $ \lambda_k = (\xi^{\alpha_{m}})^{m - k} $. To see this, let $D$ be one of the components of $ \mathcal{Y}_k $ meeting $ C $, and denote by $y$ the unique point where they intersect. Then $D$ is part of a chain of exceptional curves. Denote by $L$ the length of this chain. Using the notation and computations in Proposition \ref{prop. 3.3}, with $ C = C_{L+1} $ and $ D = C_L $, we can identify the fiber of  $ \mathcal{I}_{C} $ at $ y = y_L $ with $ <z_{L+1}> $. The eigenvalue of $ z_{L+1} $ for the automorphism induced by $\xi$ was precisely equal to $ \xi^{\alpha_{m}} $, so it follows that $ \lambda_k = (\xi^{\alpha_{m}})^{m - k} $. The result follows by re-indexing.
\end{pf}

\begin{rmk}
In particular, it is clear that this formula is independent of how we have chosen to order the elements in $ \mathcal{V} $.
\end{rmk}

\subsection{}\label{9.6}
Let us now choose an edge $ \varepsilon \in \mathcal{E} $. In the filtration of $ \mathcal{Y}_k $, we can find a subfiltration $ 0 < Z_{\varepsilon} \leq Z_{\mathcal{E}} < \mathcal{Y}_k $, with the refinements
$$ Z_{\varepsilon, L(\varepsilon) + 1} < \ldots < Z_{\varepsilon,l} < \ldots < Z_{\varepsilon,1} = Z_{\varepsilon}, $$
for any $ l \in \{ 1, \ldots, L(\varepsilon) \} $, and further refinements
$$ Z_{\varepsilon,l+1} = Z_{\varepsilon,l}^{\mu_l + 1} < \ldots < Z_{\varepsilon,l}^{k} < \ldots < Z_{\varepsilon,l}^{1} = Z_{\varepsilon,l}, $$ 
where $ Z_{\varepsilon,l}^{k} - Z_{\varepsilon,l}^{k+1} = C_{\varepsilon,l} $, for any $ k \in \{ 1, \ldots, \mu_l \} $.

As we are working with a fixed $ \varepsilon $, we will for the rest of this section suppress the index $ \varepsilon $, to simplify the notation. Take now an integer $ l \in \{1, \ldots, L-1 \} $, and let $ j_l : C_l \hookrightarrow \mathcal{Y} $ be the canonical inclusion. Consider then the subfiltration involving the component $C_l$:
$$ \ldots < Z^{\mu_l + 1}_l < \ldots < Z^k_l < \ldots < Z^1_l < \ldots . $$
At the $k$-th step in this filtration, we have $ Z_l^k - Z_l^{k+1} = C_l $ for all $ 1 \leq k \leq \mu_l $. The associated invertible sheaf at the $k$-th step is
\begin{equation} 
\mathcal{L}_l^k := j_l^*( \mathcal{I}_{Z_l^{k+1}}) = (\mathcal{I}_{C_l}|_{C_l})^{\otimes \mu_l - k} \otimes (\mathcal{I}_{C_{l+1}}|_{C_l})^{ \otimes \mu_{l+1} }. 
\end{equation}

For $ l = L $, we note that since all components in $ Z_{L} $ other than $ C_L $ have empty intersection with $ C_L $, we get instead
$$ \mathcal{L}_L^k := j_L^*( \mathcal{I}_{Z_L^{k+1}}) = (\mathcal{I}_{C_L}|_{C_L})^{\otimes \mu_L - k}. $$

Let $ g \in G $ be a group element corresponding to a primitive root of unity $ \xi $. The restriction $ g|_{C_l} $ has fixed points exactly at the two points $ y_l $ and $ y_{l-1} $ where $C_l$ meets the rest of the special fiber. We need to compute the fibers at $ y_l $ and $ y_{l-1} $ of $ \mathcal{L}_l^k $, and the corresponding eigenvalues for the automorphisms induced by $ g $ at these fibers.

Let $ g \in G $ correspond to a root of unity $ \xi \in \boldsymbol{\mu}_n $, and let 
$$ \phi_l^{k,p} : H^p(C_l, \mathcal{L}_l^k) \to H^p(C_l, \mathcal{L}_l^k) $$ 
be the automorphism induced by $ g $. We can then define the expression 
$$ Tr_l^k(\xi) := \sum_{p=0}^1(-1)^p~Tr_{\beta}(\phi_l^{k,p}). $$

\begin{ntn}
Since $ \xi^{\alpha_1} $ appears so frequently in our formulas, we introduce the notation $ \chi = \xi^{\alpha_1} $.
\end{ntn}

\begin{prop}\label{lemma 9.10}
For a primitive root $ \xi $, the expression $ Tr_l^k(\xi) $ can be computed in the following way:
\begin{enumerate}
\item If $ l \in \{2, \ldots, L - 1 \} $, we have that 
$$ Tr_l^k(\xi) = \frac{\chi^{r_{l-2}(\mu_l - k)}}{1 - \chi^{- r_{l-1}}} + \frac{\chi^{ - r_l (\mu_l - k) + r_{l-1} \mu_{l+1} }}{1 - \chi^{r_{l-1}}}, $$
for any $ k = 1, \ldots , \mu_l $.

\item If $ l = 1 $, we get
$$ Tr_1^k(\xi) = \frac{1}{1 - \chi^{ - r_0 }} + \frac{ \chi^{ - r_1 (\mu_1 - k) + r_0 \mu_2}}{1 - \chi^{r_0}} $$ 
for any $ k = 1, \ldots , \mu_1 $. 

\item Finally, if $ l = L $, we get
$$ Tr_L^k(\xi) = \frac{\chi^{ r_{L-2}(\mu_L - k)}}{1 - \chi^{ - r_{L-1}}} + \frac{1}{1 - \chi^{r_{L-1}}} $$
for any $ k = 1, \ldots , \mu_L $.
\end{enumerate}
\end{prop}
\begin{pf}
We will give the proof in case $(ii)$, when $ l \in \{2, \ldots, L - 1 \} $, and we will use the notation and results from Section \ref{lifting the action}. Recall also that $ \mathcal{L}_l^k = j_l^* \mathcal{I}_{Z_l^{k+1}} $, where 
$$ \mathcal{I}_{Z_l^{k+1}} = \mathcal{I}_{C_l}^{\otimes \mu_l - k} \otimes \mathcal{I}_{C_{l+1}}^{\otimes \mu_{l+1}} \otimes \mathcal{I}_0, $$
and where $ \mathcal{I}_0 $ has support away from $ C_l $ and $C_{l+1}$.

The fixed points of the automorphism $ g : C_l \rightarrow C_l $ are the two points $ y_{l-1}$ and $y_l$ where $C_l$ meets the other components of $ \mathcal{Z}_k $. The fibers of $ \mathcal{L}_l^k $ in the fixed points are
$$ \mathcal{L}_l^k(y_{l-1}) = \mathcal{I}_{Z_l^{k+1}}(y_{l-1}) = \mathcal{I}_{C_l}^{\otimes \mu_l - k}(y_{l-1}) = <z_{l-1}>^{\otimes \mu_l - k}, $$
and
$$ \mathcal{L}_l^k(y_l) = \mathcal{I}_{Z_l^{k+1}}(y_l) = \mathcal{I}_{C_l}^{\otimes \mu_l - k}(y_l) \otimes \mathcal{I}_{C_{l+1}}^{\otimes \mu_{l+1}}(y_l) = <w_l>^{\otimes \mu_l - k} \otimes <z_l>^{ \otimes \mu_{l+1} }. $$
Using Proposition \ref{prop. 3.3}, we compute that the eigenvalue for the automorphism on $ \mathcal{L}_l^k(y_{l-1}) $ (resp.~$ \mathcal{L}_l^k(y_l) $) is $(\chi^{r_{l-2}})^{\mu_l - k}$ (resp.~$(\chi^{ - r_l})^{\mu_l - k} (\chi^{r_{l-1}})^{\mu_{l+1}}$).

Let $ dg(y_{\star}) $ be the automorphism of the cotangent space to $C_l$ at the fixed point $y_{\star}$ induced by $g$. Using Proposition \ref{prop. 3.3} again, we compute that the eigenvalue of $ dg(y_{l-1}) $ (resp.~$ dg(y_{l}) $) is $\chi^{-r_{l-1}}$ (resp.~$\chi^{r_{l-1}}$).

We can therefore use Proposition \ref{prop. 5.6} to conclude that 
$$ Tr_l^k(\xi) = \frac{\chi^{r_{l-2}(\mu_l - k)}}{1 - \chi^{- r_{l-1}}} + \frac{\chi^{ - r_l (\mu_l - k) + r_{l-1} \mu_{l+1} }}{1 - \chi^{r_{l-1}}}. $$
\end{pf}

We can then make the following definition:

\begin{dfn}\label{contribution from edge} 
Let 
 \begin{equation} \Tr_{\varepsilon}(\xi) : = \sum_{l = 1}^{L(\varepsilon)} \sum_{k = 1}^{\mu_l} \Tr_{\varepsilon, l}^{k}(\xi). 
\end{equation}

We say that $ \Tr_{\varepsilon}(\xi) $ is the \emph{contribution to the trace} from $ \varepsilon \in \mathcal{E} $.
\end{dfn}

\begin{rmk}
Let us note that $ \Tr_{\varepsilon}(\xi) $ is defined entirely in terms of the intrinsic data of the singularity associated to $ \varepsilon $. Furthermore, this expression does not depend on the order in which we chose $ \varepsilon $. It is also clear that this expression does not depend on the chosen subfiltration of the divisor $ \sum_{i = 1}^{L(\varepsilon)} m_{\varepsilon_i} C_{\varepsilon_i} $.
\end{rmk}

\begin{rmk}
It is also easy to see that since $ \Tr_{\varepsilon}(\xi) $ is in fact a polynomial in $\xi$, the same formula is valid for any (possibly non-primitive) root of unity.
\end{rmk}

We will now show that we obtain a formula for $ \sum_{p=0}^1 (-1)^p~\Tr_{\beta}(\phi_g^p) $, where $ \phi_g^p $ is the automorphism of $ H^p(\mathcal{Y}_k, \mathcal{O}_{\mathcal{Y}_k}) $ induced by $g$, in terms of the vertex and edge contributions discussed above. 

\begin{thm}\label{thm. 9.13}
Let $ g \in G $ be a group element corresponding to a root of unity $\xi \in \boldsymbol{\mu}_n$. Then we have that
$$ \sum_{p=0}^1 (-1)^p~\Tr_{\beta}(\phi_g^p) = \sum_{\upsilon \in \mathcal{V}} \Tr_{\upsilon}(\xi) + \sum_{\varepsilon \in \mathcal{E}} \Tr_{\varepsilon}(\xi). $$ 
Furthermore, this expression depends only on the combinatorial data $ (\Gamma(\mathcal{X}_k), \mathfrak{g}, \mathfrak{m}) $ associated to $\mathcal{X}_k$. 
\end{thm}
\begin{pf}
We begin with choosing a special filtration
$$ 0 < \ldots < Z_{\varepsilon_i} < \ldots < Z_{\varepsilon_1} = Z_{\mathcal{E}} < \ldots < Z_{\upsilon_i} < \ldots < Z_{\upsilon_1} = \mathcal{Y}_k. $$
It then follows from Proposition \ref{prop-Grep} that 
$$ \sum_{p=0}^1 (-1)^p~\Tr_{\beta}(\phi_g^p) = \sum_{\upsilon \in \mathcal{V}} \Tr_{\upsilon}(\xi) + \sum_{\varepsilon \in \mathcal{E}} \Tr_{\varepsilon}(\xi), $$
where $ \Tr_{\upsilon}(\xi) $ is the expression defined in Definition \ref{contribution from vertex} and $ \Tr_{\varepsilon}(\xi) $ is the expression defined in Definition \ref{contribution from edge}.

It follows from Proposition \ref{prop. 9.7} that $ \Tr_{\upsilon}(\xi) $ only depends on the combinatorial structure of $ \mathcal{X}_k $. Likewise, Proposition \ref{lemma 9.10} gives that $ \Tr_{\varepsilon}(\xi) $ only depends on the combinatorial structure of $ \mathcal{X}_k $. Therefore, the same is true for the sum of these expressions.
\end{pf}

\subsection{An explicit trace formula}

In \cite{Thesis}, an explicit formula for $ \Tr_{\varepsilon}(\xi) $ is obtained. We present this formula here, without proof. Let us assume that $ \varepsilon $ is analytically isomorphic to the singularity $ \sigma = (m_1,m_2,n) $, with notation as in Section \ref{lifting the action}. 

Before giving the formula, we will first need to note a certain regularity of the minimal resolution of $ \sigma = (m_1,m_2,n) $, when $n$ runs through positive integers prime to $p$, and with the same residue class modulo $M$. 

\begin{prop}[\cite{Thesis}, Proposition 7.10]\label{reglocprop}
Let $m_1, m_2$ be positive integers, let $ m = \Gcd(m_1,m_2) $, and let $ M = \Lcm(m_1,m_2) $. Let us furthermore fix a positive integer $ n_0 $ that is not divisible by $p$ and that is relatively prime to $M$. Then the following properties hold:
\begin{enumerate}
\item There exists an integer $ K \gg 0 $, such that the multiplicities of the components in the minimal resolution of the singularity $ \sigma = (m_1,m_2,n) $, where $ n = n_0 + K M $, satisfy
$$ \mu_0 > \mu_1 > \ldots > \mu_{l_0} =  \ldots = m = \ldots = \mu_{L + 1 - l_1} < \ldots < \mu_{L} < \mu_{L + 1}, $$
 where $L$ denotes the \emph{length} of the singularity $ \sigma $.
\item The integers $ \mu_2 , \ldots , \mu_{l_0} $ are uniquely determined by $ \mu_0 $ and $ \mu_1 $, and similarly $ \mu_{L + 1 - l_1}, \ldots , \mu_{L-1} $ are uniquely determined by $ \mu_{L} $ and $ \mu_{L+1} $.
\item For any extension of degree $ n' = n + k M $, where $ k > 0 $, we have that the multiplicities $ \mu_l' $ of the components in the minimal resolution of the singularity $ \sigma' = (m_1,m_2,n') $ will only differ from the sequence of multiplicities associated to $ \sigma $ by inserting $m$'s ``in the middle''. 
\end{enumerate}
\end{prop}

In other words, for all $n$ sufficiently big and with a fixed residue class modulo $M$, we have that the multiplicities of the irreducible components of the exceptional components of the desingularization of $ \sigma = (m_1,m_2,n) $ are of the form as in part $ (i) $ of the proposition above. Increasing $n$ will only increase the length of the part of the components with constant multiplicity equal to $m$.

\begin{ex}\label{reglocex}
Consider the singularity $ (m_1,m_2,n) $ with $ m_1 = 3 $, $ m_2 = 4 $ and where $ n \equiv_{12} 5 $. We will use the notation $ (\mu_{L+1}, \mu_L, \mu_{L-1}, \ldots, \mu_1, \mu_0)_n $ for the multiplicities of the components in the resolution. Then we easily compute the sequences $ (3,2,3,4)_5 $, $ (3,2,1,2,3,4)_{17} $ and $ (3,2,1,1,2,3,4)_{29} $. This illustrates Proposition \ref{reglocprop} above, which then tells us that we have the sequence $ (3,2,1, \ldots, 1,2,3,4)_{5 + k \cdot 12} $, as soon as $ k \geq 1 $.
\end{ex}


In the situation where $ \varepsilon $ corresponds to the singularity $ \sigma = (m_1,m_2,n) $, where $ n \gg 0 $ with the interpretation above, it is proved in \cite{Thesis} that $ \Tr_{\varepsilon} = \Tr_{\sigma} $, where $ \Tr_{\sigma} $ is given by the formula below. 

\begin{thm}[\cite{Thesis}, Theorem 10.9]\label{Formula} 
Let $ \sigma = (m_1,m_2,n) $ be a singularity, where $ n \gg 0 $. Let $ m = \Gcd(m_1,m_2) $, and let $ \alpha_m $ (resp.~$ \alpha_{m_1} $, resp.~$ \alpha_{m_2} $) be inverse to $m$ (resp.~$m_1$, resp.~$m_2$) modulo $n$. For any root of unity $ \xi \in \boldsymbol{\mu}_n $, we have that
$$ \Tr_{\sigma}(\xi) = \sum_{r=0}^{\mu_0 - 1} (\mu_1 - \left \lceil r \frac{\mu_{1}}{\mu_{0}} \right \rceil) (\xi^{\alpha_{m_2}})^r + 
\sum_{r=0}^{\mu_{L+1} - 1} (\mu_L - \left \lceil r \frac{\mu_{L}}{\mu_{L+1}} \right \rceil) (\xi^{\alpha_{m_1}})^r - 
\sum_{r=0}^{m-1} (\xi^{\alpha_m})^r. $$
The coefficients in this expression depend only on the residue class of $n$ modulo $ \Lcm(m_1,m_2) $.
\end{thm}

With this formula at hand, one can effectively compute traces. Indeed, as long as $n$ is large compared to the multiplicities of the irreducible components of $ \mathcal{X}_k $, the expressions $ \Tr_{\varepsilon} $ can be computed using Theorem \ref{Formula}. The demand that $n$ should be ''large'' is no setback in the applications in the next section, where we are interested in the characters when $n$ grows to infinity.

\begin{ex}\label{formulaex}
Let $ m_1 = 2 $ and $ m_2 = 3 $. Then $ M = 6 $ and $ m = 1 $. So there are two cases to consider, namely $ n \equiv_6 1 $ and $ n \equiv_6 5 $. In the first case, one checks that the list of multiplicities is $ ( \mu_{L+1}, \ldots, \mu_0) = (2,1, \ldots, 1,2,3)_n $, and that $ Tr_{\sigma}(\xi) = 2 + \xi^{\alpha_3} $, where $ \alpha_3 $ is an inverse to $3$ modulo $n$.

In the second case, where $ n \equiv_6 5 $, one finds instead that the list of multiplicities is $ (2,1, \ldots, 1,3)_n $, and that the formula gives $ Tr_{\sigma}(\xi) = 1 $.
\end{ex}

\section{Character computations and jumps}\label{computations and jumps}

Let $X/K$ be a smooth, projective and geometrically irreducible curve, and let $ \mathcal{X}/S $ be the minimal SNC-model of $X$. We have in previous sections studied properties of the action of $ \boldsymbol{\mu}_n $ on the cohomology groups $ H^i(\mathcal{Y}_k, \mathcal{O}_{\mathcal{Y}_k}) $, where $ \mathcal{Y} $ is the minimal desingularization of the pullback $ \mathcal{X}_{S'} $ for some tame extension $S'/S$ of degree $n$. 

We will throughout this section make the following assumption:

\begin{ass}\label{gcdassumption}
For any $g \geq 1 $, we assume that the greatest common divisor of the multiplicities of the irreducible components of $ \mathcal{X}_k $ is $1$. If $g=1$, we assume in addition that $X/K$ has a rational point.
\end{ass}  

Let $\mathcal{J}/S$ be the N\'eron model of the Jacobian of $X$. We will in this section apply our results to the study of the filtration $ \{ \mathcal{F}^a \mathcal{J}_k \} $, where $ a \in \mathbb{Z}_{(p)} \cap [0,1] $, that we defined in Section \ref{ratfil}. We will first prove some general properties for these filtrations, and then present some computations for curves of genus $g = 1$ and $g = 2$.

We would at this point like to remark that in order to make the $ \boldsymbol{\mu}_n $-action on $ H^1(\mathcal{Y}_k, \mathcal{O}_{\mathcal{Y}_k}) $ compatible with the action on $ T_{\mathcal{J}'_k, 0} $, we have to let $ \boldsymbol{\mu}_n $ act on $ R' $ by $ [\xi](\pi') = \xi^{-1} \pi' $, for any $ \xi \in \boldsymbol{\mu}_n $. We made the choice in previous sections, when working with local rings, to let $ \boldsymbol{\mu}_n $ act by $ [\xi](\pi') = \xi \pi' $, in order to get simpler notation. This means that the irreducible characters for the representation on $ T_{\mathcal{J}'_k, 0} $ are the \emph{inverse} characters to those we compute when using our formulas for the representation on $ H^1(\mathcal{Y}_k, \mathcal{O}_{\mathcal{Y}_k}) $.

\subsection{}
Theorem \ref{thm. 9.13} states that the Brauer trace of the automorphism induced by any group element $ \xi \in \boldsymbol{\mu}_n $ on the formal difference $ H^0(\mathcal{Y}_k, \mathcal{O}_{\mathcal{Y}_k}) - H^1(\mathcal{Y}_k, \mathcal{O}_{\mathcal{Y}_k}) $ only depends on the combinatorial structure of $ \mathcal{X}_k $. If Assumption \ref{gcdassumption} is valid, we can improve this result, and get a similar result for the character of the representation of $ \boldsymbol{\mu}_n $ on $ H^1(\mathcal{Y}_k, \mathcal{O}_{\mathcal{Y}_k}) $:

\begin{thm}\label{main character theorem}
Let $ X/K $ be a smooth, projective and geometrically connected curve having genus $ g(X) > 0 $, and assume that Assumption \ref{gcdassumption} holds. Let $ \mathcal{X} $ be the minimal SNC-model of $ X $ over $S$. Furthermore, let $ S'/S $ be a tame extension of degree $n$, where $ n $ is relatively prime to the least common multiple of the multiplicities of the irreducible components of $ \mathcal{X}_k $, and let $ \mathcal{Y}/S' $ be the minimal desingularization of $ \mathcal{X}_{S'} $. 

Then the irreducible characters for the representation of $ \boldsymbol{\mu}_n $ on $ H^1(\mathcal{Y}_k, \mathcal{O}_{\mathcal{Y}_k}) $ only depend on the combinatorial data $ (\Gamma(\mathcal{X}_k), \mathfrak{g}, \mathfrak{m}) $ associated to $\mathcal{X}_k$.  
\end{thm}
\begin{pf}
Let $ g \in G $ correspond to the root $ \xi \in \boldsymbol{\mu}_n $. Then, by Theorem \ref{thm. 9.13}, we have that  
$$ \sum_{p=0}^1 (-1)^p~\Tr_{\beta}(\phi^p_g) = \sum_{\upsilon \in \mathcal{V}} \Tr_{\upsilon}(\xi) + \sum_{\varepsilon \in \mathcal{E}} \Tr_{\varepsilon}(\xi), $$
where $ \phi^p_g $ is the automorphism induced by $g$ on $ H^p(\mathcal{Y}_k, \mathcal{O}_{\mathcal{Y}_k}) $. The contributions $ \Tr_{\upsilon}(\xi) $ can be computed using Proposition \ref{prop. 9.7}, and the contributions $ \Tr_{\varepsilon}(\xi) $ can be computed by Proposition \ref{lemma 9.10}. In this way, we obtain a formula for the Brauer trace of the automorphism induced by any $ \xi \in \boldsymbol{\mu}_n $ on the formal difference $ H^0(\mathcal{Y}_k, \mathcal{O}_{\mathcal{Y}_k}) - H^1(\mathcal{Y}_k, \mathcal{O}_{\mathcal{Y}_k}) $. 

Since Assumption \ref{gcdassumption} holds, we have that $ H^0(\mathcal{Y}_k, \mathcal{O}_{\mathcal{Y}_k}) = k $ (\cite{Arwin}, Lemma 2.6). Furthermore, the $ \boldsymbol{\mu}_n $-action on $ \mathcal{Y}_k $ is relative to the ground field $k$, so it follows that the character for the representation of $ \boldsymbol{\mu}_n $ on $ H^0(\mathcal{Y}_k, \mathcal{O}_{\mathcal{Y}_k}) $ is $1$. 

We therefore obtain the formula
$$ \Tr_{\beta}(\phi^1_g) = 1 - (\sum_{\upsilon \in \mathcal{V}} \Tr_{\upsilon}(\xi) + \sum_{\varepsilon \in \mathcal{E}} \Tr_{\varepsilon}(\xi)). $$
Since the expressions $ \Tr_{\upsilon}(\xi) $ and $ \Tr_{\varepsilon}(\xi) $ only depend on the combinatorial structure of $ \mathcal{X}_k $, the same is true for $ \Tr_{\beta}(\phi^p_g) $. This completes the proof, since the Brauer character for the representation of $ \boldsymbol{\mu}_n $ on $ H^1(\mathcal{Y}_k, \mathcal{O}_{\mathcal{Y}_k}) $ is determined by the Brauer trace for the group elements $ \xi \in \boldsymbol{\mu}_n $.
\end{pf}

\vspace{0.5cm}

Let $ \mathcal{J}/S $ be the N\'eron model of the Jacobian of $X/K$. Theorem \ref{main character theorem} has the following consequence for the filtration $ \{ \mathcal{F}^a \mathcal{J}_k \} $:

\begin{cor}\label{main jump corollary}
The jumps in the filtration $ \{ \mathcal{F}^a \mathcal{J}_k \} $ with indices in $ \mathbb{Z}_{(p)} \cap [0,1] $ depend only on the combinatorial data $ (\Gamma(\mathcal{X}_k), \mathfrak{g}, \mathfrak{m}) $. In particular, the jumps do not depend on the residue characteristic $p$.
\end{cor}
\begin{pf}
Let $S'/S$ be a tame extension of degree $n$, where $n$ is prime to $l$, the least common multiple of the multiplicities of the irreducible components of $ \mathcal{X}_k $. Let $ \mathcal{J}'/S' $ be the N\'eron model of the Jacobian of $X_{K'}$. Recall from Section \ref{jacobiancase} that we could make the identification $ H^1(\mathcal{Y}_k, \mathcal{O}_{\mathcal{Y}_k}) \cong T_{\mathcal{J}'_k,0} $. 

The jumps in the filtration of $ \mathcal{J}_k $ induced by the extension $S'/S$ are determined by the irreducible characters for the representation of $ \boldsymbol{\mu}_n $ on $ T_{\mathcal{J}'_k,0} $. However, this representation is precisely the representation of $ \boldsymbol{\mu}_n $ on $ H^1(\mathcal{Y}_k, \mathcal{O}_{\mathcal{Y}_k}) $, if we let $ \boldsymbol{\mu}_n $ act on $R'$ by $ [\xi](\pi') = \xi^{-1} \pi' $, for every $\xi$. By Theorem \ref{main character theorem}, the character for this representation only depends on the combinatorial data $ (\Gamma(\mathcal{X}_k), \mathfrak{g}, \mathfrak{m}) $.  

Since $ \mathbb{Z}_{(lp)} \cap [0,1] $ is \emph{dense} in $ \mathbb{Z}_{(p)} \cap [0,1] $, we conclude that the jumps of the filtration $ \{ \mathcal{F}^a \mathcal{J}_k \} $ with indices in $ \mathbb{Z}_{(p)} \cap [0,1] $ only depend on $ \Gamma(\mathcal{X}_k) $, $ \mathfrak{g} $ and $ \mathfrak{m} $. 
\end{pf}

\vspace{0.5cm}

With the two results above at hand, we can draw some conclusions about \emph{where} the jumps occur in the case of Jacobians. Let us first recall the following terminology from \cite{Thesis}: An irreducible component $C$ of $ \mathcal{X}_k $ is called \emph{principal} if either $ p_a(C) > 0 $, or if $ C $ is smooth and rational and meets the rest of the components of $ \mathcal{X}_k $ in at least three points.  

\begin{cor}\label{specific jump corollary}
Let $ \tilde{n} $ be the least common multiple of the multiplicities of the principal components of $ \mathcal{X}_k $. Then the jumps in the filtration $ \{ \mathcal{F}^a \mathcal{J}_k \} $ occur at indices of the form $ i/\tilde{n} $, where $ 0 \leq i < \tilde{n} $. 
\end{cor}
\begin{pf}
Let us first recall that if $X$ obtains semi-stable reduction over a tame extension $K'/K$, then the Jacobian of $X$ obtains semi-abelian reduction over the same extension (see \cite{DelMum}). Furthermore, the minimal extension that gives semi-abelian reduction is the unique tame extension $ \widetilde{K}/K $ of degree $ \tilde{n} $ (\cite{Thesis}, Paper I Theorem 7.1). So in this case, the statement follows from Proposition \ref{tamejumpprop}.

Let us now assume that $X$ needs a wildly ramified extension to obtain semi-stable reduction. Consider the combinatorial data $ (\Gamma(\mathcal{X}_k), \mathfrak{g}, \mathfrak{m}) $. It follows from \cite{Winters}, Corollary 4.3, that we can find an SNC-model $ \mathcal{Z}/\Spec(\mathbb{C}[[t]]) $, where the generic fiber of $ \mathcal{Z} $ is smooth, projective and geometrically connected, and where the special fiber of $ \mathcal{Z} $ has the \emph{same} combinatorial data as $ \mathcal{X}_k $. 

Let $ \mathcal{J}_{\mathcal{Z}} $ be the N\'eron model of the Jacobian of the generic fiber of $ \mathcal{Z} $. Then the jumps of the filtration $ \{ \mathcal{F}^a \mathcal{J}_{\mathcal{Z},\mathbb{C}} \} $  occur at indices of the form $ i/\tilde{n} $, where $ 0 \leq i < \tilde{n} $. The result follows now from Corollary \ref{main jump corollary}.
\end{pf}

\subsection{}
Let $X/K$ be a smooth, projective and geometrically connected curve, and let $ \mathcal{X}/S $ be the minimal SNC-model of $X/K$. It is known that for a fixed genus $ g \geq 2 $, there are only finitely many possibilities for the combinatorial structure of the special fiber of $ \mathcal{X}/S$, modulo chains of $(-2)$-curves (\cite{Arwin},  Theorem 1.6). The same statement is, as we shall see below, also true for elliptic curves. 

Let $ \mathcal{J}/S $ be the N\'eron model of the Jacobian of $X$. Since, by Corollary \ref{main jump corollary},  the jumps of the filtration $ \{ \mathcal{F}^a \mathcal{J}_k \} $ only depend on the combinatorial structure of $ \mathcal{X}_k $, one can, for each $g > 0$, classify these jumps. In the next sections, we will give the jumps for every fiber type of genus $1$ and $2$. 

\begin{rmk}
It is not hard to see that chains of $(-2)$-curves do not affect the jumps.
\end{rmk}

\subsection{Computations of jumps for  $ g=1 $} 
Let $X/K$ be an \emph{elliptic} curve, and let $ \mathcal{E} $ be the minimal regular model of $ X $. It is a well known fact that there are only finitely many possibilities for the combinatorial structure of the special fiber $ \mathcal{E}_k $, modulo chains of $(-2)$-curves. The various possibilities were first classified in \cite{Kod}, and this is commonly referred to as the \emph{Kodaira classification}. For another treatment of this theory, we refer to \cite{Liubook}, Chapter 10.2. If now $ \mathcal{X}/S $ denotes the minimal SNC-model of $X$, it follows that there are only finitely many possibilities for the combinatorial structure of $ \mathcal{X}_k $, each one derived from the Kodaira classification. The symbols $ I, II, \ldots $ appearing in Table \ref{table 1} below are known as the \emph{Kodaira symbols} and refer to the fiber types in the Kodaira classification.

Let $ \mathcal{J}/S $ be the N\'eron model of $ J(X) = X $. It follows from Corollary \ref{main jump corollary} and Corollary \ref{specific jump corollary} that the (unique) jump in the filtration $ \{ \mathcal{F}^a \mathcal{J}_k \} $ only depends on the fiber type of $ \mathcal{X}/S $, and can only occur at finitely many \emph{rational} numbers. In Table \ref{table 1} below, we list the jumps for the various Kodaira types. Note that we obtain the same list as the one computed in \cite{Edix} by R. ~Schoof.

We would like to say a few words about how these computations are done. For each fiber type, we consider an infinite sequence $ (n_j)_{j \in \mathbb{N}} $, depending on the fiber type, where $ n_j \rightarrow \infty $ as $ j \rightarrow \infty $. For each $n_j$ in this sequence, let $ R_j/R $ be the tame extension of degree $n_j$, and let $ \pi_j $ be a uniformizing parameter of $ R_j $. Furthermore, let $ \boldsymbol{\mu}_{n_j} $ act on $ R_j $ by $ [\xi](\pi_j) = \xi \pi_j $. We can then use Theorem \ref{thm. 9.13} to compute the character for the induced representation of $ \boldsymbol{\mu}_{n_j} $ on $ H^1(\mathcal{Y}^j_k, \mathcal{O}_{\mathcal{Y}^j_k}) $, where $ \mathcal{Y}^j $ denotes the minimal desingularization of $ \mathcal{X}_{S_j} $, and where $ S_j = \Spec(R_j) $. This character is on the form $ \chi(\xi) = \xi^{i(j)} $. In particular, when $ n_j \gg 0 $, we obtain an explicit formula for $ i(j) $, using Theorem \ref{Formula}.

The character for the representation of $ \boldsymbol{\mu}_{n_j} $ on $ T_{\mathcal{J}^j_k, 0} $ is the inverse of this character, $ \chi^{-1}(\xi) = \xi^{- i(j)} $. The jump of $ \{ \mathcal{F}^a \mathcal{J}_k \} $ will then be given by the limit of the expression $ [- i(j)]_{n_j}/n_j $ as $ j \rightarrow \infty $, where $ [- i(j)]_{n_j} \equiv_{n_j} - i(j) $, and $ 0 \leq [- i(j)]_{n_j} < n_j $.

In Example \ref{example genus 1} below, we explain in detail how these computations are done for fiber type $IV$ in the Kodaira classification.

\begin{ex}\label{example genus 1}
Let $ \mathcal{X}/S $ have fibertype $IV$. In this case, the combinatorial data of $ \mathcal{X}_k $ consists of the set of vertices $ \mathcal{V} = \{ \upsilon_1, \ldots, \upsilon_4 \} $, where $ \mathfrak{m}(\upsilon_i) = 1 $ for $ i \in\{ 1,2,3 \}$, and $ \mathfrak{m}(\upsilon_4) = 3 $. Furthermore, we have that $ \mathfrak{g}(\upsilon_i) = 0 $ for all $i$. The set of edges is $ \mathcal{E} = \{ \varepsilon_1, \varepsilon_2, \varepsilon_3 \} $, where $ \varepsilon_i $ corresponds to the unique intersection point of the components $ \upsilon_i $ and $ \upsilon_4 $, for $ i = 1,2,3 $. Let us choose the ordering $ (\upsilon_i,\upsilon_4) $ for all $i$. 

Let now $ n \gg 0 $ be a positive integer relatively prime to $ p $ and to $ \Lcm( \{ \mathfrak{m}(\upsilon_i) \} ) = 3 $, and let $R'/R$ be a tame extension of degree $n$. Let $ \boldsymbol{\mu}_n $ act on $R'$ by $ [\xi](\pi') = \xi \pi' $ for any $ \xi \in \boldsymbol{\mu}_n$, where $ \pi' $ is a uniformizing parameter for $R'$. 

For any $ g \in G $, corresponding to a root of unity $ \xi \in \boldsymbol{\mu}_n $, Theorem \ref{thm. 9.13} states that
$$ \sum_{p=0}^1 (-1)^p~\Tr_{\beta}(\phi^p_g) = \sum_{\upsilon \in \mathcal{V}} \Tr_{\upsilon}(\xi) + \sum_{\varepsilon \in \mathcal{E}} \Tr_{\varepsilon}(\xi). $$

Let $ \sigma $ be the singularity $ (1,3,n) $. Then we have that $ \Tr_{\varepsilon_i}(\xi) = \Tr_{\sigma}(\xi) $ for all $ i \in \{ 1,2,3 \} $. It suffices to consider the case where $ n \equiv_3 1 $. One computes easily that $ \mu_l = 1 $ for all $ l \in \{ 1, \ldots, L(\sigma) \} $. From Theorem \ref{Formula}, we immediately get that $ \Tr_{\varepsilon_i}(\xi) = 1 $, for all $i$. 

Proposition \ref{prop. 9.7} states that 
$$ \Tr_{\upsilon}(\xi) = \sum_{k=0}^{m_{\upsilon}-1} (\xi^{\alpha_{m_{\upsilon}}})^{k} ((m_{\upsilon} - k)C_{\upsilon}^2  + 1 - p_a(C_{\upsilon})), $$
for any $ \upsilon \in \mathcal{V} $, where $ \alpha_{m_{\upsilon}} m_{\upsilon} \equiv_n 1 $. As $ C_{\upsilon_i}^2 = - 1 $ for $ i \in \{ 1,2,3 \} $, we see that $ \Tr_{\upsilon_i}(\xi) = 0 $ for these vertices, and since $ C_{\upsilon_4}^2 = - 1 $, it follows that $ \Tr_{\upsilon_4}(\xi) = - 2 - \xi^{\alpha_3} $. In total, we get
$$  \Tr_{\beta}(e(H^{\bullet}(g|_{ \mathcal{Y}_k}))) = 3 + (- 2 - \xi^{\alpha_3}) = 1 - \xi^{\alpha_3}. $$
We can therefore conclude that the character for the representation of $ \boldsymbol{\mu}_n $ on $ H^1(\mathcal{Y}_k, \mathcal{O}_{\mathcal{Y}_k}) $ is $ \chi(\xi) = \xi^{\alpha_3} $.


In order to compute the jump of the filtration $ \{ \mathcal{F}^a \mathcal{J}_k \} $, where $ \mathcal{J} $ is the N\'eron model of $J(X) = X$, we have to use the \emph{inverse} character, which is $ \chi^{-1}(\xi) = \xi^{[- \alpha_3]_n} $, where $ [- \alpha_3]_n = - \alpha_3 $ modulo $n$, and $ 0 \leq [- \alpha_3]_n < n $. The jump will be given by the limit of the expression $ ([- \alpha_3]_n)/n $ as $n$ goes to infinity over integers $n$ that are equivalent to $1$ modulo $3$. 

Since $ n = 1 + 3 \cdot h $, for some integer $h$, we get that $ \alpha_3 = \frac{1 + 2 n}{3} $, where $ 0 < \alpha_3 < n $. Therefore, the jump occurs at the limit of $ ([- \alpha_3]_n)/n = \frac{n - 1}{3 n} $ which is $ 1/3 $.

\end{ex}
 
\begin{table}[htb]\caption{Genus $1$}\label{table 1}
\begin{tabular}{|c|c|c|c|c|c|c|c|c|c|c|} 
\hline 
Fibertype & $(I)$ & $(I)^*$ & $(I_n)$ & $(I_n)^*$ & $(II)$ & $(II)^*$ & $(III)$ & $(III)^*$ & $(IV)$ & $(IV)^*$\\
\hline
Jumps & $0$ & $1/2$ & $0$ & $1/2$ & $ 1/6 $ & $ 5/6 $ & $ 1/4 $ & $ 3/4 $ & $ 1/3 $ & $ 2/3 $ \\\hline
\end{tabular}
\end{table}

\subsection{Computations of jumps for $g=2$}\label{genus 2}
Let $ X/K $ be a curve having genus equal to $2$. Like in the case for elliptic curves, there are finitely many possibilities, modulo chains of $(-2)$-curves, for the combinatorial structure of the special fiber of the minimal regular model of $X$. Moreover, there exists a complete classification of the various possible fiber types. This classification is mainly due to A.P. Ogg (\cite{Ogg}), with the exception of a few missing cases which were filled in by Y. Namikawa and K. ~Ueno in \cite{Ueno}. We use the classification and notation in \cite{Ueno}.

Let $ \mathcal{X}/S $ be the minimal SNC-model of $X$, and let $ \mathcal{J}/S $ be the N\'eron model of the Jacobian of $X$. The jumps in the filtration $ \{ \mathcal{F}^a \mathcal{J}_k \} $ depend only on the combinatorial structure of $ \mathcal{X}_k $, and can occur only at a finite set of rational numbers. 

In order to compute the jumps for each fibertype, we proceed more or less in the same manner as we did in the case of elliptic curves. In Example \ref{example 13.1}, we explain in detail how this is done for fiber type $VI$ in the classification in \cite{Ueno}. 

The jumps for the various genus $2$ fiber types are listed in tables \ref{table 2} through \ref{table 6} below. 

\begin{ex}\label{example 13.1}
We consider fiber type $VI$ in the classification in \cite{Ueno}. In this case, the set of vertices of $ \Gamma(\mathcal{X}_k) $ is $ \mathcal{V} = \{ \upsilon_1, \ldots, \upsilon_7 \} $, where $ \mathfrak{g}(\upsilon_i) = 0 $ for all $i$. Furthermore, we have that $ \mathfrak{m}(\upsilon_i) = 1 $ for $ i = 1,7 $, $ \mathfrak{m}(\upsilon_i) = 2 $ for $ i = 2,5,6 $, $ \mathfrak{m}(\upsilon_3) = 3 $ and $ \mathfrak{m}(\upsilon_4) = 4 $. The set of edges is $ \mathcal{E} = \{ \varepsilon_1, \varepsilon_2, \varepsilon_3, \varepsilon_4, \varepsilon_5, \varepsilon_6 \} $, where $ \varepsilon_1 = (\upsilon_1,\upsilon_2) $, $ \varepsilon_2 = (\upsilon_2,\upsilon_3) $, $ \varepsilon_3 = (\upsilon_3,\upsilon_4) $, $ \varepsilon_4 = (\upsilon_5,\upsilon_4) $, $ \varepsilon_5 = (\upsilon_6,\upsilon_4) $  and $ \varepsilon_6 = (\upsilon_7,\upsilon_4) $.

We have that $ \Lcm( \{\mathfrak{m}(\upsilon_i)\}) = 12 $. Let $ n \gg 0$ be any integer not divisible by $p$, and such that $ n \equiv_{12} 1 $. Let $ R'/R $ be the extension of degree $n$, and let $ \pi'$ be a uniformizing parameter of $R'$. Let $ \mathcal{Y} $ be the minimal desingularization of $ \mathcal{X}_{S'} $. We let $ \boldsymbol{\mu}_n $ act on $ R'$ by $ [\xi](\pi') = \xi \pi' $, for any $ \xi \in \boldsymbol{\mu}_n $.

Now, let $ \xi \in \boldsymbol{\mu}_n$ be a root of unity. For any $ \upsilon \in \mathcal{V} $, Proposition \ref{prop. 9.7} gives that
$$ \Tr_{\upsilon}(\xi) = \sum_{k=0}^{m_{\upsilon}-1} (\xi^{\alpha_{m_{\upsilon}}})^{k} ((m_{\upsilon} - k)C_{\upsilon}^2  + 1 - p_a(C_{\upsilon})). $$
As the computations are similar for all $ \upsilon \in \mathcal{V} $, we only do this explicitly for $ \upsilon_3 $. We have that $ p_a(C_{\upsilon_3}) = \mathfrak{g}(\upsilon_3) = 0 $, so it remains only to compute $ C_{\upsilon_3}^2 $. The edge $ \varepsilon_2 $ corresponds to the singularity $ \sigma_2 = (2,3,n) $ and the edge $ \varepsilon_3 $ corresponds to the singularity $ \sigma_3 = (3,4,n) $. Denote by $ C_l^{\sigma_2} $ the exceptional components in the resolution of $ \sigma_2 $, and by $ C_l^{\sigma_3} $ the components in the resolution of $ \sigma_3 $. Then $ C_1^{\sigma_2} $ and $ C_L^{\sigma_3} $ are the only two components of $ \mathcal{Y}_k $ that meet $ C_{\upsilon_3} $ (note the ordering of the formal branches in $ \sigma_2 $ and $ \sigma_3 $). It is easily computed that $ \mu_1^{\sigma_2} = 2 $ and that $ \mu_L^{\sigma_3} = 1 $. So it follows that $ C_{\upsilon_3}^2 = - 1 $, and therefore
$$ \Tr_{ \upsilon_3 }(\xi) = - 2 - \xi^{\alpha_3}. $$
For the other vertices, we compute that
$$ \Tr_{ \upsilon_1 }(\xi) = \Tr_{ \upsilon_7 }(\xi) = 0, $$
$$ \Tr_{ \upsilon_2 }(\xi) = \Tr_{ \upsilon_5 }(\xi) = \Tr_{ \upsilon_6 }(\xi) = - 1, $$
and
$$ \Tr_{ \upsilon_4 }(\xi) = - 7 - 5 \xi^{\alpha_4} - 3 (\xi^{\alpha_4})^2 - (\xi^{\alpha_4})^3. $$

Next, we must compute the contributions from the singularities. We will only write out the details for $ \varepsilon_3 = (\upsilon_3,\upsilon_4) $. In this case, we need to compute $ \Tr_{\sigma_3}(\xi) $. It is easily computed that $ \mu_1^{\sigma_3} = 3 $ and $ \mu_L^{\sigma_3} = 1 $. Theorem \ref{Formula} then gives that 
$$ \Tr_{\varepsilon_3}(\xi) = \Tr_{\sigma_3}(\xi) = 3 + 2 \xi^{\alpha_4} + (\xi^{\alpha_4})^2. $$
For the contributions from the other edges, we compute in a similar fashion that
$$ \Tr_{\varepsilon_1}(\xi) = \Tr_{\varepsilon_6}(\xi) = 1, $$
$$ \Tr_{\varepsilon_2}(\xi) = 2 + \xi^{\alpha_3}, $$
and
$$ \Tr_{\varepsilon_4}(\xi) = \Tr_{\varepsilon_5}(\xi) = 3 + \xi^{\alpha_4} + (\xi^{\alpha_4})^2. $$

Summing up, we get 
$$ \sum_{i=1}^7 \Tr_{\upsilon_i }(\xi) + \sum_{i=1}^6 \Tr_{\varepsilon_i}(\xi) = 1 - \xi^{\alpha_4} - (\xi^{\alpha_4})^3. $$
We can therefore conclude that the irreducible characters for the induced representation of $ \boldsymbol{\mu}_n $ on $ H^1(\mathcal{Y}_k, \mathcal{O}_{\mathcal{Y}_k}) $ are $ \chi_1(\xi) = \xi^{\alpha_4} $ and $ \chi_2(\xi) = \xi^{3 \alpha_4} $. 

The irreducible characters for the representation of $ \boldsymbol{\mu}_n $ on $ T_{\mathcal{J}'_k,0} $ induced by the action $ [\xi](\pi') = \xi^{-1} \pi' $ on $R'$ are the inverse characters of these, $ \chi_1^{-1}(\xi) = \xi^{ - \alpha_4} $ and $ \chi_2^{-1}(\xi) = \xi^{ - 3 \alpha_4} $. It is easily seen that $ [- \alpha_4]_n = (n-1)/4 $, and that $ [- 3 \alpha_4]_n = (3n-3)/4 $. Hence the jumps occur at the limits $1/4$ and $ 3/4 $ of these expressions as $ n $ goes to infinity.
\end{ex}

\begin{table}[htb]\caption{Genus $2$, Elliptic type $[1]$}\label{table 2}
\begin{tabular}{|c|c|c|c|c|c|c|}
\hline
Fiber type & $I_{0-0-0}$ & $I_{0-0-0}^*$ & $II$ & $III$ & $IV$ & $V$ \\
\hline
Jumps & $ 0 $ & $ 1/2 $ & $ 0 $, $ 1/2 $ & $ 1/3 $, $ 2/3 $ & $ 1/6 $, $ 5/6 $ & $ 1/6 $, $ 2/6 $ \\
\hline  
\end{tabular} \vspace{3 mm}

\begin{tabular}{|c|c|c|c|c|c|} 
\hline
$V^*$ & $VI$ & $VII$ & $VII^*$ & $VIII-1$ & $VIII-2$ \\
\hline
$ 4/6 $, $ 5/6 $ & $ 1/4 $, $ 3/4 $ & $ 1/8 $, $ 3/8 $ & $ 5/8 $, $ 7/8 $ & $ 1/10 $, $ 3/10 $ & $ 3/10 $, $ 9/10 $  \\
\hline  
\end{tabular} \vspace{3 mm}

\begin{tabular}{|c|c|c|c|c|c|} 
\hline
 $VIII-3$ & $VIII-4$ & $IX-1$ & $IX-2$ & $IX-3$ & $IX-4$ \\
\hline
 $ 1/10 $, $ 7/10 $ & $ 7/10 $, $ 9/10 $ & $ 1/5 $, $ 3/5 $ & $ 1/5 $, $ 2/5 $ & $ 3/5 $, $ 4/5 $ & $ 2/5 $, $ 4/5 $ \\
\hline  
\end{tabular} \vspace{3 mm}

\end{table}

\begin{table}[htb]\caption{Genus $2$, Elliptic type $[2]$}\label{table 3}
\begin{tabular}{|c|c|c|c|c|c|} 
\hline
$I_0-I_0-m$ & $I_0^*-I_0^*-m$ & $I_0-I_0^*-m$ & $2 I_0-m$ & $2I_0^*-m$ & $I_0-II-m$ \\
\hline
$ 0 $ & $ 1/2 $ & $ 0 $, $ 1/2 $ & $ 0 $, $ 1/2 $ & $ 1/4 $, $ 3/4 $ & $ 0 $, $ 1/6 $  \\
\hline  
\end{tabular} \vspace{3 mm}

\begin{tabular}{|c|c|c|c|c|} 
\hline
$I_0-II^*-m$ & $ I_0-IV-m $ & $I_0-IV^*-m$ & $I_0^*-II-m$ & $I_0^*-II^*-m$ \\
\hline
$0$, $5/6$ & $0$, $1/3$ & $0$, $2/3$ & $1/6$, $3/6$ & $3/6$, $5/6$ \\
\hline  
\end{tabular} \vspace{3 mm}

\begin{tabular}{|c|c|c|c|c|}
\hline
$I_0^*-II^*-\alpha$ & $ I_0^*-IV-m $ & $ I_0^*-IV^*-m $ & $ I_0^*-IV^*-\alpha $ & $I_0-III-m$ \\
\hline
$3/6$, $5/6$ & $1/2$, $1/3$ & $1/2$, $2/3$ & $1/2$, $2/3$ & $0$, $1/4$ \\
\hline  
\end{tabular} \vspace{3 mm}

\begin{tabular}{|c|c|c|c|c|}
\hline
$I_0-III^*-m$ & $I_0^*-III-m$ & $I_0^*-III^*-m$ & $I_0^*-III^*-\alpha$ & $ 2II-m $ \\
\hline
$0$, $3/4$ & $1/4$, $2/4$ & $2/4$, $3/4$ & $2/4$, $3/4$ & $1/12$, $7/12$ \\
\hline  
\end{tabular} \vspace{3 mm}

\begin{tabular}{|c|c|c|c|c|}
\hline
$ 2II^*-m $ & $II-II-m$ & $II-II^*-m$ & $II^*-II^*-m$ & $II^*-II^*-\alpha$ \\
\hline
$5/12$, $11/12$ & $1/6$, $1/6$ & $1/6$, $5/6$ & $5/6$, $5/6$ & $5/6$, $5/6$ \\
\hline 
\end{tabular} \vspace{3 mm}

\begin{tabular}{|c|c|c|c|c|}
\hline
$II-IV-m$ & $II-IV^*-m$ & $II^*-IV-m$ & $II^*-IV-\alpha$ & $II^*-IV^*-m$ \\
\hline
$1/6$, $2/6$ & $1/6$, $4/6$ & $2/6$, $5/6$ & $2/6$, $5/6$ & $4/6$, $5/6$ \\
\hline 
\end{tabular} \vspace{3 mm}

\begin{tabular}{|c|c|c|c|c|}
\hline
$II^*-IV^*-\alpha$ & $2IV-m$ & $2IV^*-m$ & $IV-IV-m$ & $IV-IV^*-m$ \\ 
\hline
$4/6$, $5/6$ & $1/6$, $4/6$ & $2/6$, $5/6$ & $1/3$, $1/3$ & $1/3$, $2/3$ \\
\hline
\end{tabular} \vspace{3 mm}

\begin{tabular}{|c|c|c|c|c|}
\hline
$IV^*-IV^*-m$ & $IV^*-IV^*-\alpha$ & $II-III-m$ & $II-III^*-m$ \\
\hline
$2/3$, $2/3$ & $2/3$, $2/3$ & $2/12$, $3/12$ & $2/12$, $9/12$ \\
\hline
\end{tabular} \vspace{3 mm}

\begin{tabular}{|c|c|c|c|}
\hline
$II^*-III-m$ & $II^*-III-\alpha$ & $II^*-III^*-m$ & $II^*-III^*-\alpha$  \\
\hline
$2/12$, $10/12 $ & $3/12$, $10/12 $ & $9/12$, $10/12 $ & $9/12$, $10/12 $ \\
\hline
\end{tabular} \vspace{3 mm}

\begin{tabular}{|c|c|c|c|}
\hline
$ IV-III-m $ & $ IV-III^*-m $ & $ IV-III^*-\alpha $ & $ IV^*-III-m $  \\
\hline
$3/12$, $4/12 $ & $4/12$, $9/12 $ & $4/12$, $9/12 $ & $3/12$, $8/12 $ \\
\hline
\end{tabular} \vspace{3 mm}

\begin{tabular}{|c|c|c|c|c|}
\hline
$ IV^*-III^*-m $ & $ IV^*-III^*-\alpha $ & $2III-m$ & $2III^*-m$  \\
\hline
$8/12$, $9/12 $ & $8/12$, $9/12 $ & $1/8$, $5/8$ & $3/8$, $7/8$ \\
\hline
\end{tabular} \vspace{3 mm}

\begin{tabular}{|c|c|c|c|c|}
\hline
$III-III-m$ & $III-III^*-m$ & $III^*-III^*-m$ & $III^*-III^*-\alpha$ \\
\hline
$1/4$, $1/4$ & $1/4$, $3/4$ & $3/4$, $3/4$ & $3/4$, $3/4$ \\
\hline
\end{tabular} \vspace{3 mm}

\end{table}

\begin{table}[htb]\caption{Genus $2$, Parabolic type $[3]$}\label{table 4}
\begin{tabular}{|c|c|c|c|c|c|}
\hline
$I_{n-0-0}$ & $I_n-I_0-m$ & $I_0-I_n^*-m$ & $I_n-I_0^*-m$ & $I_{n-0-0}^*$ & $I_0^*-I_n^*-m$\\
\hline
$0$ & $0$ & $ 0 $ , $ 1/2 $ & $ 0 $ , $ 1/2 $ & $ 1/2 $ , $ 1/2 $ & $ 1/2 $ , $ 1/2 $\\
\hline
\end{tabular} \vspace{3 mm}

\begin{tabular}{|c|c|c|c|c|c|}
\hline
$ II_{n-0} $ & $ II_{n-0}^* $ & $ II-I_n-m $ & $ II^*-I_n-m $ & $ IV-I_n-m $ & $ IV^*-I_n-m $ \\
\hline
$ 0 $ , $ 1/2 $ & $ 0 $ , $ 1/2 $ & $ 0 $ , $ 1/6 $ & $ 0 $ , $ 5/6 $ & $ 0 $ , $ 1/3 $ & $ 0 $ , $ 2/3 $ \\
\hline
\end{tabular} \vspace{3 mm}

\begin{tabular}{|c|c|c|c|c|c|}
\hline
$ II-I_n^*-m $ & $ II^*-I_n^*-m $ & $ II^*-I_n^*-\alpha $ & $ IV-I_n^*-m $ & $ IV^*-I_n^*-m $ \\
\hline
$ 1/6 $ , $ 3/6 $ & $ 3/6 $ , $ 5/6 $ & $ 3/6 $ , $ 5/6 $ & $ 2/6 $ , $ 3/6 $ & $ 3/6 $ , $ 4/6 $ \\
\hline
\end{tabular} \vspace{3 mm}

\begin{tabular}{|c|c|c|c|c|c|}
\hline
$ IV^*-I_n^*-\alpha $ & $ IV-II_n $ & $ IV^*-II_n $ & $ II-II_n^* $ & $ II^*-II_n^* $ & $ III-I_n-m $\\
\hline
$ 3/6 $ , $ 4/6 $ & $ 0 $ , $ 1/3 $ & $ 0 $ , $ 2/3 $ & $ 1/6 $ , $ 3/6 $ & $ 3/6 $ , $ 5/6 $ & $ 0 $ , $ 1/4 $\\
\hline
\end{tabular} \vspace{3 mm}

\begin{tabular}{|c|c|c|c|c|c|}
\hline
$ III^*-I_n-m $ & $ III-I_n^*-m $ & $ III^*-I_n^*-m $ & $ III^*-I_n^*-\alpha $ & $ III-II_n $ \\
\hline
$ 0 $ , $ 3/4 $ & $ 1/4 $ , $ 2/4 $ & $ 2/4 $ , $ 3/4 $ & $ 2/4 $ , $ 3/4 $ & $ 0 $ , $ 3/4 $ \\
\hline
\end{tabular} \vspace{3 mm}

\begin{tabular}{|c|c|c|c|c|c|}
\hline
$ III^*-II_n $ & $ III-II_n^* $ & $ III^*-II_n^* $ \\
\hline
$ 0 $ , $ 3/4 $ & $ 1/4 $ , $ 2/4 $ & $ 2/4 $ , $ 3/4 $ \\
\hline
\end{tabular} \vspace{3 mm}

\end{table}

\begin{table}[htb]\caption{Genus $2$, Parabolic type $[4]$}\label{table 5}
\begin{tabular}{|c|c|c|c|c|c|}
\hline
$I_{n-p-0}$ & $I_n-I_p-m$ & $I_{n-p-0}^*$ & $I_n^*-I_p^*-m$ & $I_n-I_p^*-m$ \\
\hline
$ 0 $ & $ 0 $ & $1/2$ & $1/2$ & $ 0 $, $1/2$ \\
\hline
\end{tabular} \vspace{3 mm}

\begin{tabular}{|c|c|c|c|c|c|}
\hline
$2I_{n}-m$ & $2I_{n}^*-m$ & $I_{n-p}$ & $III_n$ \\
\hline
$ 0 $, $1/2$ & $ 1/4 $ , $ 3/4 $ & $ 0 $, $1/2$ & $ 1/4 $ , $ 3/4 $ \\
\hline
\end{tabular} \vspace{3 mm}

\end{table}

\begin{table}[htb]\caption{Genus $2$, Parabolic type $[5]$}\label{table 6}
\begin{tabular}{|c|c|c|c|c|c|}
\hline
$I_{n-p-q}$ & $I_{n-p-q}^*$ & $II_{n-p}$ & $II_{n-p}^*$ & $III_{n}$ & $III_{n}^*$ \\
\hline
$ 0 $ & $1/2$ & $ 0 $, $1/2$ & $ 0 $, $1/2$ & $ 1/3 $, $2/3$ & $ 1/6 $, $5/6$ \\
\hline
\end{tabular} \vspace{3 mm}

\end{table}

\subsection{Final remarks and comments}
It would be interesting to know, for a curve $X/K$, the significance of the jumps in the filtration $ \{ \mathcal{F}^a \mathcal{J}_k \} $, where $ \mathcal{J} $ is the N\'eron model of $Jac(X)$. For instance, the sum of the jumps seems to be closely related to the so-called \emph{base change conductor} defined in \cite{Chai}. Furthermore, when $ k = \mathbb{C} $, and $ g = 1 $ or $ 2 $, computations show that the jumps correspond to half of the eigenvalues of the monodromy operator. It would be interesting to know if this holds in general. 

It would be nice, if possible, to have a closed formula for the irreducible characters of the representation of $ \boldsymbol{\mu}_n $ on $ H^1(\mathcal{Y}_k, \mathcal{O}_{\mathcal{Y}_k}) $, where $ \mathcal{Y} $ is the minimal desingularization of $ \mathcal{X}_{S'} $, and $ n = \Deg(S'/S) $. Such a formula would probably encode combinatorial properties of $ \mathcal{X}_k $.

We do not know if our results remain true in the case where distinct components of $ \mathcal{X}_k $ with multiplicities divisible by $p$ intersect nontrivially. The main problem is the lack of a good description of the minimal desingularization of $ \mathcal{X}_{S'} $.  

Finally, we think it would be interesting to study filtrations for N\'eron models of abelian varieties that are not Jacobians. In that case, it is not so clear what kind of data would suffice in order to determine the jumps. For instance, is it true that all jumps are rational numbers?

\bibliographystyle{plain}
\bibliography{algbib}

\end{document}